\documentclass[a4paper,12pt,twoside]{article}
\usepackage{xcolor}
\usepackage{amssymb}
\usepackage[final]{lpretex}
\usepackage{title}
\usepackage{a4}
\usepackage{amscd}
\usepackage{amsmath}
\usepackage{hyperref}
\input xy
\xyoption{all}

\EndOfDump

\def\DHpartformat#1{(#1) }
\def\DHrefpart#1{(\DHRefpart{#1})}
\LBsetstyleof{partno}{arabic}

\DocVersion[Changed to LaTeX]{1}

\let\define\def

\def\A {{\mathbb A}}  \def\C {{\mathbb C}}
  
\def\GG {{\mathbb G}}   
  
\def\N {{\mathbb N}}  \def\P {{\mathbb P}} 
\def\Q {{\mathbb Q}} \def\R {{\mathbb R}}

 \def\W {{\mathbb W}} \def\X {{\mathbb X}}
\def\Z {{\mathbb Z}} 
\def \p {\mathbb P}
\define \n {\mathbb N}
\define \z {\mathbb Z}
\define \q {\mathbb Q}
\define \PP {\mathbb P}

\def\sA {{\Cal A}} \def\sB {{\Cal B}} \def\sC {{\Cal C}}
\def\sD {{\Cal D}} \def\sE {{\Cal E}} \def\sF {{\Cal F}}
\def\sG {{\Cal G}} \def\sH {{\Cal H}} 
 \def\sK {{\Cal K}} \def\sL {{\Cal L}}
\def\sM {{\Cal M}}  \def\sO {{\Cal O}}
 \def\sR {{\Cal R}}
 \def\sT {{\Cal T}} \def\sU {{\Cal U}}
  \def\sX {{\Cal X}}

\define \cN {\Cal N}
\define \cf {\Cal F}
\define \cg {\Cal G}
\define \cE {\Cal E}
\define \ce {\Cal E}
\define \cc {\Cal C}
\define \cV {\Cal V}
\define \cA {\Cal A}
\define \cK {\Cal K}
\define \cO {\Cal O}
\define \cF {\Cal F}
\define \cn {\Cal N}
\define \cI {\Cal I}
\define \sP {\Cal P}
\def\tA {\widetilde{\Cal A}}

\def\a {\alpha} \def\b {\beta} \def\g {\gamma}  
\def\eps {\epsilon}   
\def\s {\sigma}

\define \x {\xi}
\define \y {\eta}
\define \G {\Gamma}
\define \r {\rho}
\define \w {\omega}

\def \tZ {\widetilde Z}

\def \trho {\widetilde {\rho}}

\def \tp {\widetilde{\mathbb P}}
\define \tH {\widetilde H}
\define \tG {\widetilde{G}}
\define \tW {\widetilde W}
\define \tF {\widetilde F}
\define \tm {\widetilde m}
\define \St {\widetilde S}
\define \Xt {\widetilde X}
\define \tS {\widetilde S}
\define \tpsi {\widetilde \psi}
\define \tL {\widetilde L}
\define \tE {\widetilde E}
\define \tl {\widetilde l}
\define \tA {\widetilde A}
\define \tom {\widetilde\omega}
\define \tT {\widetilde T}
\define \tB {\widetilde B}
\define \tf {\widetilde f}
\define \tsA {\widetilde{\sA}}

\define \tM {\widetilde M}
\define \tN {\widetilde N}
\define \tphi {\widetilde{\phi}}
\define \trho {\widetilde{\rho}}
\define \tR {\widetilde R}
\define \tp {\widetilde p}
\define \tq {\widetilde q}
\define \tc {\widetilde c}
\define \tsF {\widetilde {\sF}}
\define \tsP {\widetilde {\sP}}
\define \tsU {\widetilde {\sU}}
\define \th {\widetilde h}

\define \ts {\widetilde \sigma}

\def\pd {\partial}

\def \Dx1 {\frac{\pd}{{\pd} x_1}}
\def \Dy1 {\frac{\pd}{{\pd} y_1}}
\def \Dz1 {\frac{\pd}{{\pd} z_1}}
\def \Dx2 {\frac{\pd}{{\pd} x_2}}
\def \Dy2 {\frac{\pd}{{\pd} y_2}}
\def \Dz2 {\frac{\pd}{{\pd} z_2}}

\def\q {\quad} 

\def\from{\longleftarrow}

\def\mapdiagr#1{\Big\searrow\rlap{$\raise 5pt\vbox{{\hbox{$\mkern -15mu\scriptstyle#1$}}}$}}   

\def\mapdiagl#1{\llap{$\raise 5pt\vbox{{\hbox{$\scriptstyle#1\mkern
-15mu$}}}$}\Big\swarrow}              

\def\Mapdiagr#1{\nearrow\rlap{$\lower 5pt\vbox{{\hbox{$\mkern
-15mu\scriptstyle#1$}}}$}} 

\def\Mapdiagl#1{\llap{$\lower 5pt\vbox{{\hbox{$\scriptstyle#1\mkern
-15mu$}}}$}\searrow} 

\def\Mapswr#1{\swarrow\rlap{$\lower 5pt\vbox{{\hbox{$\mkern
-15mu\scriptstyle#1$}}}$}}              

\def\Mapnwl#1{\nwarrow\rlap{$\lower 5pt\vbox{{\hbox{$\mkern
-15mu\scriptstyle#1$}}}$}}

\def \inj {\hookrightarrow}

\define \Rhook {\hookrightarrow}

\def \half {\raise1pt\hbox{$\scriptstyle
        \frac{1}{2}\displaystyle$}}

\def \x{{\sl X}\llap{$\mkern -2mu {\scriptstyle -}$}}
\def \Proj {\operatorname{Proj}}
\def \Symm {\operatorname{Sym}}

\def \Pic {\operatorname{Pic}}
\def \Sing {\operatorname{Sing}}
\define \Kod {\operatorname{Kod}}
\define \dimension {\operatorname{dim}}
\define \codim {\operatorname{codim}}
\define \contr {\operatorname{contr}}
\define \rk {\operatorname{rank}}

\define \im {\operatorname{Im}}
\define \Mor {\operatorname{Mor}}
\define \Cl {\operatorname{Cl}}
\define \Hilb {\operatorname{Hilb}}

\define \degree {\operatorname{deg}}
\define \mult {\operatorname{mult}}
\define \Aut {\operatorname{Aut}}
\define \NS {\operatorname{NS}}
\define \Gal {\operatorname{Gal}}
\define \ch {\operatorname{char}}
\define \Jac {\operatorname{Jac}}
\define \Km {\operatorname{Km}}
\define \Sec {\operatorname{Sec}}
\define \Stab {\operatorname{Stab}}
\define \Br {\operatorname{Br}}
\define \inv {\operatorname{inv}}
\define \tr {\operatorname{tr}}
\define \Frob {\operatorname{Frob}}
\define \Symn {\operatorname{Sym}^n}
\define \Ev {\sE^\vee}
\define \ordp {\operatorname{ord}_p}
\define \Supp {\operatorname{Supp}}
\define \Ann {\operatorname{Ann}}
\define \disc {\operatorname{disc}}
\define \Lie {\operatorname{Lie}}
\define \embdim {\operatorname{embdim}}

\def \Fix{\operatorname{Fix}}
\def\Isom{\operatorname{Isom}}

\def\lng{\operatorname{long}}

\def\W{\cansymb{W}}

\def\ad{\operatorname{ad}}

\def\vol{\operatorname{vol}}
\def\wt{\operatorname{wt}}

\def\bZ{\overline Z}

\def\bC{\overline{C}}

\def\bu{\bf{u}}
\def\bv{\bf{v}}
\def\bw{\bf{w}}
\def\bx{\bf{x}}

\def\tsX{\widetilde{\sX}}

\begin{document}
\title
[Flag varieties and reductive groups]
{Recognizing flag varieties and reductive groups}
\author{I. Grojnowski}
\address{DPMMS\\
  CMS Wilberforce Rd.\\
  Cambridge  CB3 0WB U.K.}
\email{groj@dpmms.cam.ac.uk}

\author{N. I. Shepherd-Barron}
\address{Dept. of Mathematics\\
King's College\\
Strand\\
London WC2R 2LS U.K.}
\email{Nicholas.Shepherd-Barron@kcl.ac.uk}
\maketitle
\begin{abstract} Fix a flat and projective morphism $X\to\Sigma$
  of schemes. We show, first, that any set of $\P^1$-fibrations on $X$
  defines a set of simple roots,
  a set of simple coroots and a Cartan matrix
  $C$. Second, $X$ is an {\'e}tale
  $\sF$-bundle over some projective $\Sigma$-scheme,
  where $\sF$ is the flag variety of the
  adjoint Chevalley group over $\Sp\Z$ defined by $C$.
  In particular, if the simple roots generate $\NS(X/\Sigma)_\Q$
  and $X$ is cohomologically flat in degree zero over $\Sigma$ then
  $X$ is a form of $\sF$. When $X$ is a smooth Fano variety
  over $\Sp\C$
  all of whose extremal rays are accounted for by these fibrations
  this is due to Occhetta, Sol{\'a}~Conde, Watanabe and Wi{\'s}niewski.
  Third, we
  recover, in a uniform way, the isomorphism and isogeny theorems
  of Chevalley and Demazure: over any base a pinned reductive
  group is determined by its pinned root datum, and a $p$-morphism
  of pinned root data determines a unique homomorphism
  of the corresponding groups.
\end{abstract}
\begin{section}{Introduction}
  We assume throughout this paper that
  $\Sigma$ is a connected scheme and that
  $X\to\Sigma$ is flat and projective.
  
  If there is a set of $\P^1$-fibrations
  $\{\pi_i:X\to Y_i : i\in I\}$ where each $Y_i$ is also defined
  over $\Sigma$ then we say that
  $X$ is \emph{multiply fibred}. Our main result is this:

  \begin{theorem}\label{main} If $X$ is multiply fibred then
    $X$ is an {\'e}tale $\sF$-bundle over some flat projective
    $\Sigma$-scheme,
    where $\sF$ is the flag variety of the adjoint Chevalley group over
    $\Z$ that is determined by certain intersection numbers on $X$.
    \noproof
  \end{theorem}
  When $\Sigma=\Sp\C$, $X$ is a smooth Fano variety and the
  given $\P^1$-fibrations account for all the extremal rays of $X$,
  so that $X=\sF$,
  this is due to
  \cite{OSCWW}. Their argument involves some delicate
  geometry
  in the cases that lead to $F_4$ and $G_2$; ours builds upon
  their ideas but is uniform. It works by getting
  better control over the cohomology of line bundles on
  the Bott--Samelson varieties $Z$ constructed from $X$
  by regarding $Z$
  as a member of a family that is parametrized by $X$.
  In fact, our approach explains that Bott--Samelson varieties
  and their collapsing maps to the flag variety are combinatorial
  objects defined purely in terms of the Cartan matrix, despite
  the presence of moduli in their construction
  as iterated marked $\P^1$-bundles.
  
  To begin, we show that a multiply conic ind-projective
  variety $X$ leads to a generalized Cartan matrix $C$.
  When $X$ is projective and the conic bundle structures
  are all smooth $C$ turns out to be a Cartan matrix.
  We then construct Bott--Samelson varieties and use them to
  prove the main result. As a consequence we recover, without
  any individual consideration
  of groups of semi-simple rank $\le 2$,
  the theorems of Chevalley and Demazure \cite{SGA3}:
  a pinned reductive group over any base is determined by its
  pinned root datum and a $p$-morphism of pinned root data
  determines a unique homomorphism between the corresponding
  groups.
\end{section}
\begin{section}{Multiple conic bundle structures}\label{Cartan}
  We suppose that
  $X$ is ind-projective and \emph{multiply conic} in that there is a set
  $\{\pi_i:X\to Y_i : i\in I\}$ of conic bundle structures.
  That is, $\pi_i$ is projective and flat and its
  geometric generic fibres are copies of $\P^1$.
  We let $\a_{X, i}=\a_i$ denote
  the dual of the relative dualizing sheaf of $\pi_i:X\to Y_i$
  and $\a^\vee_{X, i}=\a^\vee_i$ the class of a geometric fibre of $\pi_i$.
  We see next that the $\a_i$ and $\a^\vee_j$ can be regarded
  as the simple roots and simple coroots of some Kac--Moody
  group.

  \begin{proposition} \label{cart}
    The intersection matrix $C=(C_{i j})=(\a_i.\a^\vee_j)$
    is a generalized Cartan matrix. That is, $C_{i i}=2$, $C_{i j}\le 0$
    if $i\ne j$ and $C_{i j}=0$ if and only if $C_{j i}=0$.
    \begin{proof} For this we can suppose
      that $\Sigma$ is a geometric point.

      $C_{i i}=2$ is clear. So suppose $i\ne j$.
      Choose a geometric generic point $y\in Y_i$
      and consider $\pi_j^{-1}(\pi_j(\pi_i^{-1}(y)))$.
      This is a surface whose image under $\pi_j$ is a curve in $Y_j$.
      Pulling back by the normalization $B$ of this curve
      gives a surface $\bZ$, a morphism $g:\bZ\to B$
      and a curve $\bC$ in $\bZ$ which is the fibre $\pi_i^{-1}(y)$.
      The map $g$ is induced
      by $\pi_j$: its fibres are fibres of $\pi_j$.
      So the geometric generic fibre of $g$
      is $\P^1$ and every geometric fibre is reduced with only nodes.
      \footnote{If $X$ is multiply fibred rather than merely multiply
        conic then $\bZ\to B$ is a $\P^1$-bundle. So $\bZ$ is smooth
        and the argument that follows is even easier.}
      So $\bZ$ has only RDPs, of type $A$. The curve
      $\bC$ is a multi-section of $g$ and is collapsed by
      the non-trivial morphism $\pi_i$.

      Let $Z\to \bZ$ be the minimal resolution
      and $C$ the strict transform of $\bC$. Since
      $C$ is collapsed by $\pi_i$, $C^2\le 0$.

      Since $B$ is rational, we have
      $K_C=K_B+R$
      where $\deg R\ge 0$, even if
      $C\to B$ is inseparable. Also
      $K_C=(K_Z+C)\vert_C$ and
      $K_Z=K_{Z/B}+K_B =-\a_j+K_B$.
      So
      $$0\le \deg R=-\a_j.C+C^2,$$
      so that
      $(\a_j. \a_i^\vee)=\a_j.C\le C^2\le 0$.

       If $(\a_j. \a_i^\vee)= 0$ then the formation
      of $\bZ$ is unchanged by switching $i$ and $j$,
      so $\a_i.\a^\vee_j=0$.
    \end{proof}
  \end{proposition}
  
    Let $s_i$ be the reflexion in $\NS(X)_\R$ defined by
    the pair $(\a_i. \a_i^\vee)$
    and $W=\langle S\rangle$.

  \begin{proposition} $(W, S)$ is a Coxeter system.
    \begin{proof} This is a consequence of
      Proposition \ref{cart} and the results of \cite{Go}.
    \end{proof}
  \end{proposition}
  We shall see later (Theorem \ref{3.7}) that when $X$ is
  finite-dimensional  and multiply fibred  then
  the system $(W, S)$ is finite, so that $C$ is positive definite.
  That is, $C$ is a Cartan matrix.
\end{section}
\begin{section}{Fibrations and Weyl groups}\label{six}
  Suppose that $X$ is multiply fibred.
  We know that the set of $\P^1$-fibrations
  gives rise to a generalized Cartan matrix $C$
  and a Coxeter system $(W, S)$. We shall prove here
  that $(W, S)$ is finite and crystallographic and deduce
  that $C$ is a Cartan matrix.
  For this we can assume that $\Sigma$ is a geometric point.
  We follow \cite{OSCWW} very closely.

  Say $\dim X=m$ and put
  $\omega_{X/Y_i}=\sO_X(K_i)=-\a_{X, i}=-\a_i$.

  \begin{lemma}\label{2.4+}
    For any
    $D\in\Pic(X)$ the two sheaves
    $\pi_{i*}\sO(D+(D.\a^\vee_i+1)K_{i}))$ and
    $(\pi_{i*}\sO(K_i-D))^\vee$
    on $Y_i$ are locally free and have equal Chern polynomials.
    \begin{proof} The local freeness is clear.
      Put $D.\a^\vee_i=l$. Both sheaves are zero if
      $l\ge -1$, so assume $l\le -2$. 
       Suppose first that
      $\pi_i$ is a Zariski bundle, so that
      there is a Cartier divisor $H$ on $X$ such that
        $H.\a^\vee_i=1$. Then $\pi_{i*}\sO(H)=\sF$, say,
        is locally free
        of rank $2$ and $X=\P(\sF)$.
        Then $D\sim lH+\pi_i^*B$, some $B\in \Pic({Y_i})$,
        and $K_i\sim -2H+\pi_i^*\det\sF$. Then
        $K_i-D\sim -2H+\pi_i^*\det\sF-lH-\pi_i^*B$,
        so that
        $$\pi_{i*}\sO(K_i-D)\cong\Symm^{-(l+2)}\sF
        \otimes B^\vee\otimes\det\sF.$$
        Also
        $D+(l+1)K_i\sim lH+\pi_i^*B-
        2(l+1)H+\pi_i^*\det\sF^{l+1},$
        so that
        $$\pi_{i*}\sO(D+(l+1)K_i)\cong\Symm^{-(l+2)}\sF
        \otimes B\otimes\det\sF^{\otimes(l+1)}.$$
        It is enough to observe that
        $$c_t((\Symm^{-(l+2)}\sF)^\vee)=
        c_t(\Symm^{-(l+2)}\sF\otimes\det\sF^{\otimes(l+2)}),$$
        which is a consequence of the
        splitting principle. 
        So the result is proved if $\pi_i$ is a Zariski bundle.

        Quite generally, the projection formula shows that two locally free
        sheaves $\sA,\sB$ on a projective $k$-scheme $Y_i$
        have equal Chern polynomials if there
        is a proper morphism $f:Z\to Y_i$ such that
        $f_*\sO_Z=\sO_{Y_i}$ and $c_t(f^*\sA)=c_t(f^*\sB)$.
        Take $f:Z\to Y_i$ to be $\pi_i:X\to Y_i$; then
        $pr_1:X\times_{Y_i} X\to X$ has a section,
        so is a Zariski $\P^1$-bundle. Then the sheaves
        in question have equal Chern classes after pulling back under $pr_1$
        and we are done.
    \end{proof}
  \end{lemma}

  \begin{corollary}\label{2.4++}
    $\chi(Y_i, (\pi_{i*}\sO(K_i-D))^\vee)=
    \chi(Y_i, \pi_{i*}\sO(D+(D.\a^\vee_i+1)K_i)).$
    \begin{proof} This follows from
      Lemma \ref{2.4+}
      and
      the splitting principle.
    \end{proof}
  \end{corollary}

  We omit the proofs of the remaining results in this section
  since they are exactly
  as in \cite{OSCWW}.
  
  We let $N^1=N^1(X)$ denote $\Pic(X)$ modulo numerical
  equivalence and $N_1=N_1(X)$ denote the dual
  group of curves modulo numerical equivalence.
  The Euler characteristic is then a polynomial function
  $\chi=\chi_X:N^1\to \Z$ of degree $m$.
  Suppose that $\rho_i\in N^1_\Q$ such that
  $\rho_i.\a^\vee_i=1$.
  Consider the shift
  $$T_i:N^1_\Q\to N^1_\Q: D\mapsto D+\rho_i$$
  and the function $\chi^{T_i}=\chi^{T_i}_X=\chi\circ T_i^{-1},$
  which is also polynomial of degree $m$. There
  are simple reflexions
  $$r_i:N^1 \to N^1:D\mapsto D+(D.\a^\vee_i)K_i.$$
  Put $S=\{r_i\vert i\in I\}$ and
  $W=\langle S\rangle \subseteq GL(N^1)$.
  We know that $(W, S)$ is a Coxeter system.
  Since $N^1$ is a $\Z$-lattice, in order to show that
  $(W, S)$ is finite and crystallographic it is enough
  to prove that $W$ is finite.
  
  For $w\in W$ let $w^{T_i}$ denote
  the affine linear transformation defined by
  $$w^{T_i}=T_i^{-1}\circ w\circ T_i.$$
  
  \begin{lemma}\label{2.7}\label{2.3+}\label{1.3}
    $\chi^{T_i}(r_i(D))=-\chi^{T_i}(D)$.
    \noproof
  \end{lemma}

  Consider the volume polynomial
  $\vol=\vol_X:N^1\to \Z$
  defined by $\vol(D)=c_1(D)^m$. This is
  a homogeneous polynomial of degree $m$
  and is, up to a factor of $1/m!$,
  the leading part of each of the functions
  $\chi$ and $\chi^{T_i}$.

  \begin{corollary}\label{homog}\label{1.7}
    $\vol_X\circ w=\det (w)\vol_X$
    for all $w\in W$.
    \noproof
  \end{corollary}

Each reflexion $s$ acting on $N^1_\R$
has a mirror $H_s$, its fixed locus $\Fix(s)$,
which is defined by a homogeneous linear form $L_s$
that is unique up to a scalar.
If $s=r_i$ then $H_s=(\a^\vee_i)^\perp.$
\begin{lemma}\label{dim} There are at most $m$ mirrors in $N^1_\R$.
  \begin{proof} The linear form $L_s$ belonging to a mirror $H_s$
    divides the polynomial $\vol_X$.
  \end{proof}
\end{lemma}
For $s\in S$, corresponding to $\a_s^\vee$, write
$D(s)=\{x\in N^1_\R : (x.\a_s^\vee)>0\}$. For $I\subseteq S$
put $\sC_I=\cap _{s\in I}D(s)$ and
$\sD_I=\cap_{s\in I}\Fix(s)\cap\sC_{S-I}.$ Put
$T=\cup_{w\in W}wSw^{-1}$.
\begin{lemma}\label{distinct}
  Distinct elements of $T$
  have distinct fixed loci.
  \begin{proof} It is enough to show that if $s, t\in S$,
    $q$ is a conjugate of $t$ and $\Fix(s)\subseteq \Fix(q)$
    then $q=s$.

    Put $I=\{s\}$. Then $\sD_I$ is not empty and $q, s$
    act trivially on $\sD_I$. In particular
    $q(\sD_I)\cap s(D_I)\ne\emptyset$, so that,
    by \cite{Bo}, p. 96, Prop. 5, $qW_I=sW_I$.
    Then $1\ne q\in W_I=\{1, s\}$, so
    that $q=s$.
  \end{proof}
\end{lemma}
\begin{theorem}\label{3.7}
  \part[i] $(W, S)$ is a finite
  crystallographic Coxeter system.
  \part[ii] The generalized Cartan matrix
  $C$ is positive definite.
  \begin{proof} There are only finitely many mirrors
    in $N^1_\R$, and so, by Lemma \ref{distinct},
    $T$ is finite. According to \cite{Bo},
    p. 14, Lemme 2, $\ell(w)\le \vert T\vert$, so that $\ell$
    is a bounded function and $W$ is finite.
    Since $N^1$ is a $\Z$-lattice $W$ is also crystallographic.

    The positivity of $C$ is a well known consequence.
  \end{proof}
\end{theorem}
Of course, a generalized Cartan matrix that is positive definite is a
Cartan matrix.
\end{section}
\begin{section}{Bott--Samelson varieties}
  \begin{subsection}{Preliminaries}\label{SectPrelim}
    Let $X\to\Sigma$ be multiply fibred; this
    determines a Cartan matrix $C$, and a finite
    Weyl group $W$, as in Section \ref{six}.

    Now let $G=G(C)$ the corresponding split
  adjoint Chevalley group over $\Z$;
  we take the existence of $G$ for granted.
  We let $\sF=\sF_\Z$ denote
  the flag scheme of $G$ (the scheme of
  Borel subgroups); this
  has $\P^1$-fibrations $\tau_i:\sF\to\sP_i$.
  We identify the simple roots $\a_i$ with the
  relative tangent bundles of the various $\tau_i$
  and the simple coroots $\a_i^\vee$
  with the classes of their fibres.
  As for further notation, $B, T, U, W$ have their usual meaning;
  $Q=\X^*(T)\subseteq\Pic(\sF)$, the root lattice; $\Pi=\{\a_i\}$,
  which is a root basis; $\Phi=W(\Pi)$ is the set of roots;
  $\Phi_\pm$ is the set of positive (negative) roots;
  $\Phi_{++}=\Phi_+-\Pi$. 
  We normalize things by requiring the weights of
  $\Lie(U)$ to be the \emph{negative} roots. That is, is $U$ generated
  by the copies $U_{-\a}$ of the additive group $\GG_a$ with $\a\in \Pi$.

  In particular, $\sF\to\Sp\Z$ is multiply fibred. Moreover,
  $C$, $W$ and the simple roots $\a_{\sF,i}$ defined in
  Sections \ref{Cartan} and \ref{six} for $\sF$ coincide
  with the objects $C$, $W$ and the simple roots $\a_i$
  mentioned in the previous paragraph,
  and also any $\phi\in Q$ defines a line bundle on $\sF$
  and a line bundle on $X$, both denoted by $\phi$.

  Suppose that ${\bw}=(s_{i_1},\ldots ,s_{i_n})$
  is a word of
  simple reflexions and $w=s_{i_1}\ldots s_{i_n}$.
  We say that $\bw$ \emph{represents} $w$. If
  $\ell(w)=n$ then $\bw$ is \emph{reduced}.
  The Bott--Samelson
  variety associated to $X$ and $\bw$ is
  $$\tZ^X_{\bw}=X\times_{Y_{i_1}}\times X\times
  \cdots \times_{Y_{i_n}}X.$$
  The first projection $q^X_{\bw}:\tZ^X_{\bw}\to X$
  is an iterated $\P^1$-bundle. Our convention is that
  if ${\bw}$ is empty then $\ell({\bw})=0$
  and $q^X_{\bw}$ is the identity morphism.
  We denote by $f^X_{\bw}:\tZ^X_{\bw}\to X$ the last projection.

  If $s_\a$
  is a simple reflexion and ${\bw}={\bv}s_\a$ is a word,
  then there is a $\P^1$-bundle
  $p^X_{\bw}:\tZ^X_{\bw}\to\tZ^X_{\bv}$ and a section
  $i^X_{\bw}:\tZ^X_{\bv}\inj \tZ^X_{\bw}$ of $p^X_{\bw}$
  that fit into a commutative diagram
  $$
  \xymatrix
  {
    &{\tZ^X_{\bw}}\ar[r]^{f^X_{\bw}}\ar[d]^{p^X_{\bw}}\ar@/_/[dl]_{q^X_{\bw}}&
    {X}\ar[d]^{\pi_\a}\\
    {X}&{\tZ^X_{\bv}}\ar[l]^{q^X_{\bv}}\ar@/^/[u]^{i^X_{\bw}}
    \ar[ur]_{f^X_{\bv}}\ar[r]_g&
    {Y_\a}
  }
  $$
  where the square is Cartesian.
  Concatenation of words defines a composition 
  $$\tZ^X_{\bx}\times_{f^X_{\bx}, X, q^X_{\bw}}\tZ^X_{\bw}\to
  \tZ^X_{{\bx}{\bw}}$$
  so that, if we define $\tZ^X_{\W}=\lim_{\to}\tZ^X_{\bw}$,
  then $(q, f):\tZ^X_{\W}\to X\times_\Sigma X$
  is a semi-groupoid in the category of strict
  ind-projective $\Sigma$-schemes. For $x\in X$
  we let $X(x)\subseteq X$ denote the equivalence class
  of $x$ defined by this ind-semi-groupoid. 
  We shall return to
  this later.

  A $\P^1$-bundle is \emph{marked} if it has a distinguished section.
  Giving such a bundle on a scheme $Y$ is equivalent to giving
  a line bundle $\sL$ on $Y$ and a class $\xi\in H^1(Y, \sL^{-1})$
  defined up to units: the marked bundle is $\P(\sE)$ where
  $\sE$ is the extension of $\sO_Y$ by $\sL^{-1}$ defined
  by $\xi$. The section is defined by the surjection
  $\sE\to\sO_Y$.
  
  Suppose that ${\bv}=(s_{\a_1},\ldots ,s_{\a_n})$
  is a word. We say that \emph{$s_\a$ occurs in $\bv$}
  if $\a$ is one of the $\a_i$.
  
  If $\phi\in Q$, so is a line bundle on $X$ or $\sF$,
  then we also let $\phi$ denote the line bundle $f^{X *}_{\bw}\phi$
  on $\tZ^X_{\bw}$ and $f^{\sF *}_{\bw}\phi$ on
  $\tZ^\sF_{\bw}$.

  We denote by $\sK^X_{\bw}$ the $X$-group scheme
  of automorphisms of $q^X_{\bw}:\tZ^X_{\bw}\to X$ that
  preserve the iterated marked $\p^1$-bundle structure of
  $\tZ^X_{\bw}$. 

  \begin{lemma} Every line bundle on $\tZ^X_{\bw}$
    is $\sK^X_{\bw}$-linearized.
    \begin{proof} Induction on $\ell({\bw})$. If $\ell({\bw})=0$
      there is nothing to prove, so suppose ${\bw}={\bv}s_\a$.
      The action of $\sK^X_{\bw}$ on the given section of the
      marked $\P^1$-bundle $p^X_{\bw}:\tZ^X_{\bw}\to\tZ^X_{\bv}$
      defines a homomorphism $\sK^X_{\bw}\to \sK^X_{\bv}$,
      so line bundles on $\tZ^X_{\bw}$
      that are pulled back from $\tZ^X_{\bv}$
      are $\sK^X_{\bw}$-linearized. The fact that the given
      section is preserved by $\sK^X_{\bw}$ shows that
      the relative $\sO(1)$ is $\sK^X_{\bw}$-linearized
      and the lemma is proved.
    \end{proof}
  \end{lemma}
  The marked $\P^1$-bundle
  $(p^X_{{\bv}s_\a}, i^X_{{\bv}s_\a})$
  is defined by a class
  $\xi^X\in H^1(\tZ^X_{\bv}, -\a)$
  modulo units.
  Given any morphism $x\to X$ with fibre
  $\tZ^X_{\bv}(x)= (q^X_{\bv})^{-1}(x)$ we can restrict
  $\xi^X$ to $H^1(\tZ^X_{\bv}(x), -\a)$
  and ask whether this restriction is zero.

  \begin{lemma}\label{basic}
    $\xi^X$ is nowhere zero on $X$ if 
    $s_\a$ occurs in $\bv$.
    \begin{proof} Suppose that $s_\a$ occurs in $\bv$.
      Say $s_\a=s_{i_t}$. Then
      $Z^{X}_{\{s_{i_t},s_\a\}}=\P^1\times\P^1
      \inj Z^{X}_{{\bv}s_\a}$, and hence the restriction of the
      $\P^1$-bundle $Z^{X}_{{\bv}s_a}\to Z^{X}_{\bv}$ to
      the subscheme $Z^{X}_{s_{i_t}}$ of $Z^{X}_{\bx}$
      is the trivial $\P^1$-bundle. Then the
      restriction of $\xi^{X}$ to $H^1(Z^{X}_{s_{i_t}},-\a)$
      is non-zero, since it is represented by an extension
      $$0\to \sO\to\sO(1)\oplus\sO(1)\to\sO(2)\to 0$$
      on $Z^{X}_{s_{i_t}}=\P^1$, so $\xi^X$ is nowhere zero.
    \end{proof}
  \end{lemma}
\end{subsection}
\medskip

{\bf{For the rest of this section we fix a word $\bv$ and make the
    following
    induction hypothesis:
    for each word $\bx$ that is obtained from $\bv$ by a
    truncation from the right (including the vacuous
    truncation ${\bx}={\bv}$)
    the morphism $q^X_{\bx}:\tZ^X_{\bx}\to X$ is an
    {\'e}tale $Z^\sF_{\bx}$-bundle, where
    $Z^\sF_{\bx}=(q^\sF_{\bx})^{-1}(z)$
    for some $z\in\sF(\Z)$.}}

  More precisely, the inductive hypothesis is that
  there is a smooth affine cover $\sX=\sX_{\bv}=\Sp\sO_\sX\to X$,
  a morphism $\sX\to\sF$ and an $\sX$-isomorphism
  $\tZ^X_{\bv}\times_X\sX \to \tZ^\sF_{\bv}\times_\sF\sX$.
  (Note that $\tZ^\sF_{\bv}\to\sF$ is such a bundle
  since it can be identified with $Z^\sF_{\bv}(z)\times^B G$.)

  Note that this holds when $\bv$ is the empty word.
\begin{subsection}{$Q$-filtered complexes}\label{S = Fk}
  Recall that a filtration $F$ of a complex $M$
  is a sequence of complexes
  $$
  \xymatrix
  {
    {\cdots}\ar[r]&
    {F^{i+1}}\ar[r]&
    {F^{i}}\ar[r]&
    {\cdots}\ar[r]&
    {M.}    
  }
  $$
  It is \emph{separated} if $\lim_{\from} F^i=0$
  and \emph{exhaustive} if $\lim_{\to} F^i\to M$ is a
  quasi-isomorphism. Furthermore it is \emph{finite}
  if $F^n\to M$ is a quasi-isomorphism for all $n\le N_-$
  and $F^n=0$ for all $n\ge N_+$ for some $N_-\le N_+$ in
  $\Z$.
  We set $gr_F^i(M)=cone(F^{i+1}\to F^{i})$
  and $gr_F(M)=\oplus gr^i(M)$.
  A finite filtration 
  induces a convergent spectral sequence,
  which we refer to as \emph{the spectral sequence of
    the filtered complex $M$}, whose $E_1$ page is given by
  $$E^{i j}_1= \sH^{i+j}(gr^i(M))\Rightarrow \sH^{i+j}(M).$$
  If $F\to M$ is a finite filtration and for every $i\in\Z$
  there is given a finite filtration ${}^iG$ of
  $gr_F^i(M)$ then there is a well defined
  finite filtration $\tF\to M$, the \emph{refinement}
  of $F$ by $\{ {}^iG\}$, such that
  $$gr_{\tF}^j(M)= gr_{{}^iG}^{j-N_i}(gr^i_F(M)),$$
  where $N_i < j \le N_{i+1}$ for some increasing
  sequence of integers $N_i$. To see this,
  suppose given ${}^iG^j\to cone (F^{i+1}\to F^i)$. Then define
  $$\tG^j=cone ({}^iG^j\oplus F^i\to cone(F^{i+1}\to F^i))[1],$$
  giving $F^{i+1}\to\tG^j\to F^i$. Iterating this constructs
  $\tF$.

  If $f:X\to Y$ is a proper morphism and $F\to M$ is a finite
  filtration of a complex of coherent sheaves on $X$
  then $Rf_*F\to Rf_*M$ is a finite filtration and
  $gr^i Rf_*M\cong Rf_*gr^i M$.

  A filtration $F$ of a complex $M$ of sheaves on $\tZ^X_{\bv}$
  is a \emph{$Q$-filtration} if it is finite, $\sK^X_{\bv}$-linearized
  and each $gr_F^i(M)$ is isomorphic to $\g_i$
  for some $\g_i\in Q$. We call these $\g_i$ the \emph{weights}
  of $M$ and the set of these weights, counted with multiplicity,
  is denoted $\wt(gr(M))$. The cohomology sheaves of $M$
  might have torsion, however.
\end{subsection}
\begin{subsection}{The key lemma}
  Recall that, if $\b,\delta$ are roots, the $\delta$-string
  through $\b$, denoted by $I_\delta(\b)=
  (\Phi\cup\{0\})\cap (\b+\Z\delta)$,
  is an interval that contains 
  $\b$ and $s_\delta(\b)$.
\begin{lemma}
  (the key lemma) Suppose that $\bx$ is a word
  obtained by truncating $\bv$ from the right.
  \begin{enumerate}
  \item If $M$ is a $Q$-filtered complex on $\tZ^X_{\bx}$
    then $Rq^X_{{\bx}*}M$ is a $Q$-filtered
    complex on $X$.
  \item If $\b\in\Phi_{++}$ then
    $\wt(gr(Rq^X_{{\bx}*}(-\b)))\subseteq -\Phi_{++}$.
    \item If $\a\in\Pi$ then 
      $\wt(gr(Rq^X_{{\bx}*}(-\a)))\subseteq -\Phi_{++}\cup\{-\g_0\}$
      where $\g_0=0$ if $s_\a$ appears in ${\bx}$ and $\g_0=\a$ otherwise.
      Moreover, 
      $-\g_0$ has multiplicity $1$.
    \end{enumerate}
  \begin{proof}
      We argue by induction on $\ell({\bx})$. If
      $\ell({\bx})=0$ then $\tZ^{X}_{\bx}={X}$ and
      the lemma is immediate.      
      Otherwise suppose that ${\bx}={\bu}s_\delta$
      and that the lemma
      has been proved for $\bu$. Then
      consider the marked $\P^1$-bundle
      $p^X_{\bx} : \tZ^{X}_{\bx} \to \tZ^{X}_{\bu}$.
      
      Suppose $\b\in\Phi_+\cup\{0\}$ and consider
      $R{p^X_{\bx}}_*(-\b)=\sE$.

      \begin{enumerate}
      \item $\b=0$: then $\sE=\sH^0(\sE)=\sO$, the zero weight.

      \item $\b=\delta$: then $\sE=\sH^1(\sE)=\sO$.

      \item $\b\in \Phi_+\setminus\{\delta\}$.
        Say $(\b.\delta^\vee)=-n$.
        \begin{enumerate}
        \item $n\ge 0$. Then $\sE=\sH^0(\sE)$ and is filtered;
          the set $gr(\sE)$ is
          $$\{-\b, -\b-\delta, \ldots , -\b-n\delta\}
          \subseteq [s_\delta(-\b), -\b]\subseteq I_\delta(-\beta).$$
          This interval is contained in $\Phi\cup\{0\}$,
          and we claim that these roots lie in $-\Phi_{++}$.

          Since $\b\ne\delta$, there exists simple $\g\in\Supp(\b)$
          with $\g\ne\delta$. So $\g, \delta\in \Supp(\b+i\delta)$
          for $i\in [1, n]$, and so every $\b+i\delta\in\Phi_{++}$
          and the claim is established.

          So every piece of $gr(\sE)$ lies in $-\Phi_{++}\cup\{-\b\}$
          in this case.
        \item $n=-1$. Then $\sE=0$ and there no weights.
        \item $n\le -2$. Since $(\b.\delta^\vee)\le 0$,
          $\b\in\Phi_{++}$.
          This time $\sE=\sH^1(\sE)$ and $gr(\sE)$ is
          $$\{-\b+\delta, \ldots , -\b-(n+1)\delta\} \subsetneq
          [-\b, s_\delta(-\b)]\subseteq I_\delta(-\b),$$
          so again they are roots (none of them is zero,
          clearly). Since $\b$ is not simple there exists
          a simple root $\g\in\Supp(\b)$ with $\g\ne \delta$.
          Then $\g,\delta\in\Supp(\b+i\delta)$
          for $-1\ge i\ge n+1$.
        \end{enumerate}
      \end{enumerate}
      So 
      $gr(Rp^X_{{\bx}*}(-\b))$ consists of $0$ if
      $\b=0$ or $\delta$ and 
      elements of $-\Phi_{++}\cup\{-\b\}$ otherwise.
      More generally, for any $\g\in Q$ it is clear that
      $Rp^X_{{\bx}*}(-\g)$ is a $Q$-filtered complex, concentrated
      in one degree.

      We know that $R{{p^X_{{\bx}*}}}M$ has a filtration whose graded
      pieces are ${{Rp^X_{{\bx}*}}}(-\g)$
      for various $-\g\in \wt(gr(M))$. Refining this filtration
      via the above $Q$-filtration of $R{{p^X_{{\bx}*}}}(-\g)$
      endows $R{{p^X_{{\bx}*}}}M$ with a $Q$-filtration.
      Since $Rq^X_{{\bx}*}M=R{{p^X_{{\bx}*}}}Rq^X_{{\bu}*}M$,
      the induction hypothesis shows that
      $Rq^X_{{\bx}*}M$ is $Q$-filtered
      and the lemma is proved.
    \end{proof}
  \end{lemma}

  \begin{lemma}\label{**} If $\a\in\Pi$ then
    each non-zero term
    ${}^XE^{i, j}_\infty =gr^i(R^{i+j}q^X_{{\bv}*}(-\a))$
    in the spectral sequence of the filtered complex
    $Rq^X_{{\bv}*}(-\a)$
    is of the form $\sM_\b(-\b)$
    where $\b\in \Phi_{++}\cup\{0, \a\}$,
    $\sM_\b$ is the pullback of a coherent sheaf on $\Sigma$
    and $\sM_0=\sO$.
    \begin{proof} First consider the case where $X=\sF$.
      On $\sF$
      the sheaf $\sO(-\a)$ is $G$-linearized
      and the $Q$-filtered complex $Rq^{\sF}_{{\bv}*}(-\a)$
      is $G$-linearized. However, for distinct $\lambda, \mu\in Q$
      and any finite complexes of sheaves $\sM_\lambda,\sM_\mu$
      of sheaves that are pulled back from $\Sp\Z$,
      $$Hom_{\sF}(\sM_\lambda\otimes\lambda,
      \sM_\mu\otimes\mu)^G=0.$$
      This forces
      the differentials in the spectral sequence
      to preserve $Q$-weights, as this is equivalent to requiring that,
      if
      ${}^\sF E^{i, j}_r=\sM'(\g),$ then
      $d^{i j}_r({}^\sF E^{i, j}_r) =0$
      unless
      $${}^\sF E^{i+r, j-r+1}_r =\sM_\g''(\g)$$
      for some coherent sheaves $\sM_\g', \sM_\g''$
      pulled back from $\Sp\Z$. In particular
      ${}^\sF E^{i, j}_\infty$ is of the required form.

      On $X$ we argue as follows. By the key lemma
      each ${}^XE^{i j}_1$ is a line bundle
      pulled back from $X$
      and has the same $Q$-weights as ${}^\sF E^{i, j}_1$.
      By the induction assumption, there is an isomorphism
      between the pullbacks to $\sX$ of
      $R^iq^X_{{\bv}*}(-\a)$
      and $R^iq^\sF_{{\bv}*}(-\a)$.
      
      It follows that each differential in
      each page
      ${}^X E^{. .}_r$ preserves $Q$-weights,
      as whether a
      differential ${}^Xd^{i, j}_r$ is non-zero
      can be detected {\'e}tale locally on $X$
      and the ${}^\sF d^{i, j}_r$
      do preserve weights, as just proved.
      Since the $Q$-weights of ${}^\sF E^{i, j}_1$
      are the same as those of ${}^X E^{i, j}_1$
      this remains true for each ${}^\sF E^{i, j}_r$
      and ${}^X E^{i, j}_r$.
    \end{proof}
  \end{lemma}

  \begin{corollary}\label{cor} Suppose that $\a\in\Pi$.
    
    \part[i] $H^0(X, R^iq^X_{{\bv}*}(-\a))=0$ for all $i\ne 1$.

    \part[ii] If $s_\a$ occurs in $\bv$ then $H^0(X,
    R^1q^X_{{\bv}*}(-\a))$
    is free of rank $1$ and is generated by the class $\xi^X$.

    \part[iii] If $s_\a$ does not occur in $\bv$ then $H^0(X,
    R^1q^X_{{\bv}*}(-\a))=0$.
    \begin{proof} If $\b\in\Phi_+$ there is a simple root $\delta$
      such that $(\b.\delta^\vee)>0$ and hence
      $R^0\pi_{\delta  *}(-\b)=0$, where $\pi_\delta:X\to Y_\delta$
      is the $\P^1$-fibration corresponding to $\delta$,      
      so $H^0(X,-\b)=0$ for all $\b\in\Phi_+$.
      Then Lemma \ref{**} shows that in the $Q$-filtered complex
      $M=Rq^X_{{\bv}*}(-\a)$ there is at most one $gr^i(M)$
      (the piece of weight zero, if it exists) that has non-zero
      global sections, and if it does it 
      contributes to $H^0(X, R^1 q^X_{{\bv}*}(-\a))$.
      The corollary is proved once we remark that, by
      Lemma \ref{basic}, the class $\xi^X$ generates
      $H^0(X, R^1 q^X_{{\bv}*}(-\a))$, whether this module is
      zero or not.
    \end{proof}
  \end{corollary}
\end{subsection}
\begin{subsection}{Comparing $X$ with $\sF$}
  We fix $X\to\Sigma$, multiply fibred.
  We now complete the induction to prove that
  for any word $\bw$
  the morphism $q^X_{\bw}:\tZ^X_{\bw}\to X$ is an {\'e}tale
  $Z^\sF_{\bw}$-bundle, where $Z^\sF_{\bw}=
  (q^\sF_{\bw})^{-1}(z)$ for any $\Z$-point
  $z$ of $\sF$. 
  In the notation of the induction hypothesis,
  we identify the $\sX$-schemes
  $\tZ^X_{\bv}\times_X\sX$, $\tZ^\sF_{\bv}\times_\sF\sX$
  and $Z^\sF_{\bv}\times_{\Sp\Z}\sX$. Similarly
  we identify the $\sX$-fibres of $\sK^X_{\bv}$ and $\sK^{\sF}_{\bv}$
  with the constant $\sX$-group $K^\sF_{\bv}\times_{\Sp\Z}\sX$
  where $K^\sF_{\bv}=\Aut(Z^\sF_{\bv})$.
  
  \begin{lemma}\label{restrict}
    If
    $f:Z\to X$ is a morphism of schemes,
    $\sG$ is a quasi-coherent sheaf
    on $Z$ and $\sX$ is an affine scheme with a morphism
    $\sX\to X$
    then the restriction map
    $H^1(Z,\sG)\to H^1(Z\times_X \sX,\sG)$
    factors through
    $H^0(X, R^1f_*\sG)$.
    \begin{proof}  The Leray spectral sequence
      gives      a commutative     diagram
      $$\xymatrix{
        {H^1(X, f_*\sG)}\ar[r]\ar[d] &
        {H^1(Z, \sG)}\ar[r]\ar[d] &
        {H^0(X, R^1f_*\sG)}\ar[d]^\delta\ar[r]&
        {H^2(X, f_*\sG)}\ar[d]\\
        {H^1(\sX, f_*\sG)}\ar[r]&
        {H^1(Z\times_X \sX,\sG)}\ar[r]^-\epsilon&
        {H^0(\sX, R^1f_*\sG)}\ar[r]&
        {H^2(\sX, f_*\sG)}
      }$$
      whose rows are exact. Since $H^i(\sX,f_*\sG)=0$ for all
      $i\ge 1$ the map
      $\epsilon$ is an isomorphism.
      The lemma follows, via $\epsilon^{-1}\circ\delta$.
    \end{proof}
  \end{lemma}

  Let $X'=X$ or $\sF$
  and consider the restriction map
  $r^{X'}:H^1(\tZ^{X'}_{\bv}, -\a)\to H^1(Z^{X'}_{\bv}\times\sX, -\a)$.

  \begin{proposition}\label{prelim}
    Assume that $H^0(X', R^1 q^{X'}_{{\bv}*}(-\a))\ne 0.$

    \part[i] $R^1q^{X'}_{{\bv}*}(-\a)^{\sK^{X'}_{\bv}}\cong\sO_{X'}.$

    \part[ii] $H^0(X', R^1q^{X'}_{{\bv}*}(-\a))=
    H^0(X', (R^1q^{X'}_{{\bv}*}(-\a))^{\sK^{X'}_{\bv}})$.

    \part[iii] $r^X(\xi^X)$ and $r^{\sF}(\xi^{\sF})$
    generate the same $\sO_{\sX}$-module in
    $H^1(Z^\sF_{\bv}\times\sX, -\a)$.
    \begin{proof}
      \DHrefpart{i} Suppose $\g\in Q$ and
      $\g^{\sK^\sF_{\bv}}$ is not 
      the zero sheaf. Then $\g^T=\g$ where $T$
      is a maximal torus in $B$ and we consider
      the action of $B$ on $\tZ^X_{\bv}$ via the
      natural homomorphism
      $\rho^X_{\bv}:B\to\sK^X_{\bv}$. This induces
      a homomorphism $\pi^X_{\bv}:T\to \sK^X_{\bv}/\sU^X_{\bv}$
      where $\sU^X_{\bv}$ is the unipotent radical of
      $\sK^X_{\bv}$. The kernel $\ker\pi^X_{\bv}$ is generated
      by the cocharacters $\b^\vee$ for which $s_\b$ does
      not occur in $\bv$. So if $\g^{\sK^X_{\bv}}\ne 0$
      then $(\g.\b^\vee)=0$ for all $s_\b$ occurring in $\bv$.

      By assumption $s_\a$ occurs in $\bv$. Then,
      by Lemma \ref{**}, only the term $\sO$
      can contribute to the space of $\sK^{X'}_{\bv}$-invariant
      global sections of $R^1q^X_{{\bv}*}(-\a)$,
      as for all other terms $\sM_\g(-\g)$ there is an
      $s_\b$ occurring in $\bv$ for which
      $(\g.\b^\vee)\ne 0$.

      \DHrefpart{ii} This is an immediate consequence.

      \DHrefpart{iii} Both
      $r^X(\xi^X)$ and $r^{\sF}(\xi^{\sF})$
      lie in the $\sO_\sX$-module
      $(H^1(Z_{\bv}\times\sX,-\a))^{K_{\bv}\times\sX}$.
  \end{proof}
\end{proposition}
\begin{proposition}\label{bundle}
  For any word $\bw$ the morphism
  $q^X_{\bw}:\tZ^X_{\bw}\to X$ is an {\'e}tale
  $Z^{\sF}_{\bw}$-bundle.
  \begin{proof} We can suppose that ${\bw}={\bv}s_\a$.
    If $s_\a$ does not occur in $\bv$ then
    $\tZ^X_{\bw}=\Proj_{\tZ^X_{\bv}}(\sO\oplus\sO(\a))$
    and there
    is nothing to prove. Otherwise we argue as follows.

    Both      $r^X(\xi^X)$ and $r^{\sF}(\xi^{\sF})$
    lie in the $\sO_\sX$-module
    $(H^1(Z_{\bv}\times\sX,-\a))^{K_{\bv}\times\sX}$.
    This $\sO_\sX$-module is free of rank one and,
    by Lemma \ref{basic},
    each is a generator of it.
    We take $\sX_{{\bw}}$ to be the
    projectivization of the line bundle
    over $\sX_{\bv}$ that is generated
    by $\xi^X$, so that
    $\sX_{\bw}=\sX_{\bv}$.
  \end{proof}
\end{proposition}
\end{subsection}  
\end{section}
\begin{section}{Geometry}\label{geometry}
  We suppose that $\Sigma$ is connected,
  that $X$ is projective and flat
  over $\Sigma$
  and that $X$ is multiply fibred. We write $\sF_\Z=\sF$.
  
  Given a word $\bw$, not necessarily reduced,
  there is a commutative diagram
  $$
  \xymatrix
  {
    {\tZ^X_{\bw}}\ar[r]^{p_{\bw}}\ar@/_1pc/[drr]_{q^X_{\bw}}\ar@/^1.5pc/[rr]^{r_{\bw}}&
    {V_{\bw}}\ar[r]^{\nu_{\bw}}\ar[dr]_{\s_{\bw}}&
    {R_{\bw}}\ar[d]^{\rho_{\bw}}\ar@{^{(}->}[r]&
    {X\times_\Sigma X}\ar[dl]^{pr_2}\\
    &&{X}&
  }
  $$
  where the top row is the Stein factorization of
  $(f^X_{\bw}, q^X_{\bw}):\tZ^X_{\bw}\to X\times X$.
  In particular, $\nu_{\bw}$ is finite and the homomorphism
  $\sO \to\nu_{{\bw}*}\sO$ is injective; that is,
  $\nu_{\bw}$ is finite and dominant.
  We know, by Proposition \ref{bundle},
  that $q^X_{\bw}$
  is an {\'e}tale $Z^\sF_{\bw}$-bundle. That is, there
  is an {\'e}tale cover $\sX\to X$ and an $\sX$-isomorphism
  $$h:Z^\sF_{\bw}\times_{\Sp\Z}\sX
  \buildrel\cong\over\to
  \tZ^X_{\bw}\times_X\sX.$$
  
  Recall that
  the image $\sF_w$ in $\sF$ of $Z^\sF_{\bw}$
  under $f^\sF_{\bw}:\tZ^\sF_{\bw}\to\sF$
  depends only on the element $w$ of $W$
  that is represented by $\bw$ and is normal
  relative to $\Sp\Z$
  \cite{An},
  \cite{RR}, \cite{Se}.
  Let $E^\sF_{\bw}\subseteq Z^\sF_{\bw}$ denote the exceptional
  locus of the restriction of $f^\sF_{\bw}$ to $Z^\sF_{\bw}$ and
  $\tE^X_{\bw}\subseteq\tZ^X_{\bw}$ the exceptional
  locus of $f^X_{\bw}:\tZ^X_{\bw}\to X$.
  Then the isomorphism $h$ takes $E^\sF_{\bw}\times_{\Sp\Z}\sX$
  to $\tE^X_{\bw}\times_X\sX$
  since each of these is the locus covered by those curves that
  are orthogonal to the group generated by the line bundles that
  are represented by
  the simple roots. Observe also that the exceptional locus
  of $Z^\sF_{\bw}\times_{\Sp\Z}\sX\to R_{\bw}\times_X\sX$
  is $E^\sF_{\bw}\times_{\Sp\Z}\sX$.

  \begin{proposition}\label{5.1} $\s_{\bw}:V_{\bw}\to X$ is an {\'e}tale
    $\sF_w$-bundle.
    \begin{proof} There is
      a commutative diagram
      $$
      \xymatrix
      {
        {Z^\sF_{\bw}\times_{\Sp\Z}\sX}\ar[r]^-h\ar[drr]&
        {\tZ^X_{\bw}\times_X\sX}\ar[r]^-{r_{\bw}}\ar[dr]&
        {R_{\bw}\times_X\sX}\ar[d]\\
        &&{\sX}
      }
      $$
      We have just seen that the exceptional locus of
      $Z^\sF_{\bw}\times_{\Sp\Z}\sX\to R_{\bw}\times_X\sX$
      is $E^\sF_{\bw}\times_{\Sp\Z}\sX$.
      Write $\tZ^X_{\bw}\times_X\sX=\tsX$.
      Then $r_{\bw}\circ h$ factors through
      $\sF_w\times_{\Sp\Z}\sX$ and so $r_{\bw}\vert_{\tsX}$ does too.
      Since $p_{\bw}\vert_{\tsX}$ is characterized
      as the minimal proper morphism of $\sX$-schemes
      such that (i) its source is $\tsX$, (ii)
      $(p_{\bw}\vert_{\tsX})_*\sO=\sO$
      (this is true here because $\sF_w$ is normal)
      and (iii) $p_{\bw}\vert_{\tsX}$ collapses exactly those curves
      that are collapsed by $r_{\bw}\vert_{\tsX}$,
      it follows that $p_{\bw}\vert_{\tsX}$
      factors through $\sF_w\times_{\Sp\Z}\sX$ also,
      and then that $V_{\bw}\times_X\sX$
      is $\sX$-isomorphic to $\sF_w\times_{\Sp\Z}\sX$.
    \end{proof}
  \end{proposition}

  Suppose that
  ${\bw}_0$ represents the longest element $w_0$ of $W$
  and define $R\subseteq X\times_\Sigma X$ to be the equivalence relation
  generated by the ind-semi-groupoid
  $(f, q):\tZ^X_{\W}\to X\times_\Sigma X$.

  \begin{lemma} $R_{{\bw}_0} =R$.
    \begin{proof} If $x=\Sp k(x)\in X$
      is any field-valued point then
      $pr_1^{-1}(x)\cap R$ is the equivalence class $[x]$,
      which is irreducible. Also $\dim [x]\le \ell(w_0)$,
      since also $[x]$ is dominated by $\cup \sF_w$,
      so that $[x]=\rho_{{\bw}_0}^{-1}(x)$.
    \end{proof}
  \end{lemma}

  So we have a commutative diagram
  $$
  \xymatrix
  {
    {V_{{\bw}_0}}\ar[dr]_{\s}\ar[r]^\nu&{R}\ar[d]^\rho\\
    &{X}
  }
  $$
  where $\nu=\nu_{{\bw}_0}$,
  $\rho=\rho_{{\bw}_0}$ and
  $\s={\s_{{\bw}_0}}$ is an {\'e}tale $\sF$-bundle.

  For any $x=\Sp k(x)\in X$ let $X(x)=[x]_{red}$.

  \begin{lemma} $X(x)$ is geometrically integral over $k(x)$.
    \begin{proof} It is dominated by the smooth and irreducible
      $k(x)$-scheme
      $Z^X_{{\bw}_0}\otimes k(x)$.
    \end{proof}
  \end{lemma}

  \begin{corollary} There is a finite and dominant morphism
    $\psi:\sF\otimes k(x)\to X(x)$.
    \noproof
  \end{corollary}

  \begin{proposition} For each $i$ there is a Cartesian diagram
    $$
    \xymatrix
    {
      {\sF\otimes k(x)}\ar[r]^-\psi\ar[d]_{\tau_i\otimes 1_{k(x)}}&
      {X(x)}\ar@{^(->}[r]\ar[d]^{\pi_i(x)}&
      {X}\ar[d]^{\pi_i}\\
      {\sP_i\otimes k(x)}\ar[r]&
      {Y_i(x)}\ar[r]&
      {Y_i.}
    }
    $$
    That is, the $\P^1$-fibrations on $X$ restrict to
    $\P^1$-fibrations on $X(x)$.
    \begin{proof} There are morphisms
      $$
      \xymatrix
      {
        {\sF\otimes k(x)}\ar[r]^-{\psi}&
        {X(x)}\ar@{^(->}[r]&
        {X}\ar[r]^{\pi_i}&
        {Y_i.}
      }
      $$
      Define $\pi_i(x):X(x)\to Y_i(x)$ to be
      the Stein factorization of $X(x)\to Y_i$.      Observe that
      $\tau_i\otimes 1_{k(x)}:\sF\otimes k(x)\to\sP_i\otimes k(x)$
      is then the Stein factorization of $\sF\otimes k(x)\to Y_i(x)$,
      so that we have constructed the required
      commutative diagram.

      Since $\tau_i\otimes 1_{k(x)}$ and $\pi_i$ are $\P^1$-fibrations
      and the geometric fibres of $\tau_i\otimes 1_{k(x)}$
      map isomorphically to geometric
      fibres of $\pi_i$ the outer rectangle
      is Cartesian.

      For the right hand square, consider the Zariski tangent
      space $T(z)$ at a point $z$ of $X(x)$. The subspace of $T(z)$
      that is annihilated
      under the composite map $X(x)\to X\to Y_i$
      has dimension at most $1$. But $\pi_i(x)$
      has $1$-dimensional fibres and so
      kills at least a $1$-dimensional subspace of $T(z)$.
      Therefore $\pi_i(x)$ is smooth, and therefore the
      right hand square is Cartesian.

      It follows that the left hand square is also Cartesian.
    \end{proof}
  \end{proposition}
  \begin{lemma}\label{mult}
    
    \part[i] Suppose that, for each $i$,
    $V_i$ is a proper closed algebraic subset
    of $\sP_i$ and that $\tau_i^{-1}(V_i)=\tau_j^{-1}(V_j)$
    for all $i, j$. Then every $V_i$ is empty.
    \part[ii] $X(x)$ is smooth over $k(x)$.
    \begin{proof} \DHrefpart{i} Say $V=\tau_i^{-1}(V_i)$.
      Then the class $[V]$ in $H^*=H^*(\sF, \Q_{\ell})$
      is invariant
      under each simple reflexion, and so under $W$.
      But $H^*$ is the regular representation of $W$,
      so that $[V]$ is a multiple of $[\sF]$
      and then must vanish.

      \DHrefpart{ii} Take $V_i=\b_i^{-1}(\Sing Y_i(x))$
      and use \DHrefpart{i}.
    \end{proof}
  \end{lemma}

  \begin{theorem}\label{finite} 
    $\psi:\sF\otimes k(x)\to X(x)$ is
    an isomorphism.
    \begin{proof} Abbreviate $\sF\otimes k(x)$ to $\sF$.
      We know that $\psi$ is
      finite and dominant and
      induces an isomorphism
      $\psi^*:N^1(X(x))_\Q\to N^1(\sF)_\Q$ that takes
      simple roots to simple roots. Moreover, $\psi$
      maps simple coroots to simple coroots.
      Put $\rho=c_1(X(x))/2$ so that,
      in the notation of the discussion following
      Corollary \ref{2.4++}, we can take $\rho_i=\rho$
      and write $\chi^{T_i}_{X(x)}=\chi^T_{X(x)}.$
      Given a coroot $\a^\vee$, define the
      linear form $F_{\a^\vee}$ on $N^1$
      by $F_{\a^\vee}(D)=(D.\a^\vee)/(\rho.\a^\vee).$
      Then by Lemma \ref{1.3}, $\chi^T_{X(x)}$
      is divisible, as a polynomial over $\Q$, by
      $\prod_{\a^\vee>0}F_{\a^\vee}$. Since both polynomials
      have the same degree (namely, $\dim\sF$),
      it follows that
      $\chi^T_{X(x)}=\lambda \prod_{\a^\vee>0}F_{\a^\vee}$
      for some $\lambda\in\Q$. Since
      $\chi_{X(x)}(\sO_{X(x)})=\chi_{X(x)}^T(\rho)=\lambda,$
      it follows that $\lambda\in\Z$.
      Similarly $\chi_\sF^T=\mu \prod_{\a^\vee>0}F_{\a^\vee},$
      and since $\chi_\sF(\sO_\sF)=1$ we get $\mu=1$.
      (Of course, this is nothing but the Weyl dimension
      formula.)
      So $\chi_{X(x)}^T/\chi_\sF^T\in\Z$.
      So $\vol_{X(x)}/\vol_\sF\in\Z$.
      But $\vol_\sF/\vol_{X(x)}=\deg\psi$,
      so that $\deg\psi=1$.
      Since $X(x)$ is smooth, the theorem is proved.
    \end{proof}
  \end{theorem}

  \begin{theorem}\label{5.8} $\nu:V_{{\bw}_0}\to R$
    is an isomorphism.
    \begin{proof} We know that for any field-valued point $x$
      of $X$ the morphism
      $\s^{-1}(x)\to\rho^{-1}(x)_{red}$
      is an isomorphism. Therefore $\nu$ separates points.
      Similarly, it separates tangent vectors.
      Since it is proper, it is 
      a closed embedding. Since it is dominant
      it is an isomorphism.
    \end{proof}
  \end{theorem}
  
  \begin{theorem}\label{9.1}
    
    \part[i] $X$ is an {\'e}tale $\sF$-bundle over
    some projective $\Sigma$-scheme $H$.

    \part[ii] Suppose that the simple roots $\a_i$ span
    $\NS(X/\Sigma)_\Q$ and that $X\to\Sigma$ is
    cohomologically flat in degree zero.
    Then $X$ is an {\'e}tale $\sF$-bundle over 
    $\Sigma$.
    \begin{proof}\DHrefpart{i} By Proposition \ref{5.1}
      and Theorem \ref{5.8} the subscheme
      $R$ of $X\times X$ is an {\'e}tale $\sF$-bundle over $X$.
      In particular, it is flat and projective over $X$
      and so defines an $X$-point of $\Hilb_{X/\Sigma}$.
      Let $H$ denote the scheme-theoretic image
      of the classifying morphism $X\to \Hilb_X$;
      then $X\to H$ is the quotient $X\to X/R$.
      Therefore $X\to H$ is smooth and the
      fibre product $X\times_H X\to X$ is
      an {\'e}tale $\sF$-bundle. Therefore $X\to H$
      is an {\'e}tale $\sF$-bundle.

      \DHrefpart{ii} The hypotheses imply that $\NS(H/\Sigma)$
      is finite, so that $H$ is finite over $\Sigma$. The
      cohomological flatness of $X\to\Sigma$ in degree zero
      gives $H=\Sigma$.
    \end{proof}
  \end{theorem}
\end{section}
\begin{section}{Uniqueness for pinned reductive groups}
  We recover Chevalley
  and Demazure's theorem concerning the uniqueness of pinned
  reductive groups. In the next section we shall recover their
  more general homomorphism theorem (or isogeny theorem) as a corollary.

  We remark that
  the existence over $\Sp\Z$ of such a group for a given root datum
  is proved in \cite{SGA3} as a consequence of their uniqueness;
  however, we shall assume existence over $\Sp\Z$
  since there are now proofs of this that are
  independent of the uniqueness theorem.

  Suppose that $G$ is a reductive group over $\Sigma$,
  $\sF^G$ its flag scheme and $C(G)$ its Cartan matrix,
  which is determined by intersection numbers
  on $\sF^G$.

  \begin{lemma}\label{adj} Suppose that $H$ is an algebraic group
    (that is, a smooth connected affine
    group scheme of finite type) over a field $k$ that acts
    transitively on some non-empty projective $k$-variety $X$
    and that $K$ is the kernel of the action. 

    \part[i] If $K^0$ is a subgroup of a torus
    then $H$ is reductive.

    \part[ii] If $K^0$
    is a finite subgroup of a torus then $H$ is semi-simple.

    \part[iii] If $K=1$ and $\chi(X, \sO_X)=1$ then $H$
    is adjoint.
    \begin{proof} The hypotheses and
      conclusions are insensitive to field extension
      so we can assume that $k$ is algebraically closed.

      Let $S$ denote the radical of $H$.
      Then the fixed locus $V$ of the $S$-action is not empty.
      Since $S$ is a normal subgroup $V$ is preserved by $H$
      and so $V=X$. That is, $S\subseteq K^0$ and so is a torus,
      so that $H$ is reductive.

      If $K^0$ is also finite then $S=1$ and $H$ is
      semi-simple.

      Finally, suppose that $K=1$ and $\chi(X, \sO_X)=1$ .
      Then the centre $Z$ of $H$ is finite and
      acts freely on $X$, so that
      $1=\chi(X, \sO_X)=\deg(Z)\chi(X/Z, \sO_{X/Z}).$
    \end{proof}
  \end{lemma}
  Recall that $\chi(X, \sO_X)=1$ when $X$ is a flag variety.

  \begin{theorem}\label{rigid}\cite{De}
    $\sF^G$ is rigid and $\Aut^0(\sF^G)=G^{ad}$.
    \begin{proof}  We can assume that $C(G)$ is irreducible,
      that $G$ is adjoint and that $\Sigma$ is
      the spectrum of an algebraically closed field $k$.

      Put $\sF^G=\sF$. Suppose that $\a$ is a simple root
      and that $\tau_\a:\sF\to \sP_\a$ is the corresponding
      $\P^1$-bundle. This gives homomorphisms
      $$H^1(\sF, \Theta_\sF)\to
      H^1(\sF, \tau_\a^*\Theta_{\sP_\a})=
      H^1(\sP_\a, \Theta_{\sP_\a})$$
      from which it follows that every deformation $\tsF$ of $\sF$ over
      $\Sp k[\eps]$ contracts to a deformation $\tsP_\a$ of $\sP_\a$ over
      $\Sp k[\eps]$. Since $\tsF$ and $\tsP_\a$ are
      locally isomorphic to $\sF\otimes\Sp k[\eps]$
      and $\sP_\a\otimes\Sp k[\eps]$, respectively,
      the morphism $\tsF\to\tsP_\a$ is smooth.
      Therefore $\tsF$ is multiply fibred and is then trivial,
      by Theorem \ref{9.1}. So the tangent space
      $H^1(\sF,\Theta_\sF)$ to moduli vanishes and
      the rigidity is proved.

      This $H^1$ is also the obstruction space to the
      smoothness of $\Aut^0({\sF})=L$, say, so that
      $L$ is smooth. 
      There is an ample $L$-linearized line bundle on $\sF$
      (for example, $\omega_{\sF}^{-1}$)
      and so $L$ is affine. Certainly, $L$
      acts effectively on $\sF$; since $L$
      contains $G$ it also acts transitively.
      Therefore $L$ is adjoint.

      So $\sF=G/B$ and $\sF=L/Q$
      where $Q$ is a parabolic subgroup scheme of $L$. Consider
      the projection $\pi:L/Q_{red}\to L/Q$; this is
      $L$-equivariant, and so, since $G\subseteq L$,
      is $G$-equivariant. Because $\pi$ is radicial
      $L/Q_{red}$ is homogeneous under $G$,
      and so $L/Q_{red} =G/P$ for some parabolic subgroup scheme
      $P\subseteq B$. So $P=B$ and 
      $Q=Q_{red}$.

      Decomposing $C(L)$ into its connected
      components corresponds to breaking $L$
      into simple factors $L_i$.
      Then $G/B = \prod L_r/Q_r$
      and we see that $L$ is simple. Let $\Delta$ denote a root basis
      for $L$; then $Q$ corresponds to a proper
      subset $I$ of $\Delta$ which is empty if and only if
      $Q$ is a Borel subgroup of $L$. 

      Suppose that $\Phi_\Delta$ is the root system generated by $\Delta$
      and $\Phi_I\subset \Phi_\Delta$ that generated by $I$.
      Then, for all $\delta\in \Delta-I$,
      $$\Phi_\Delta\cap (\Phi_I+\Q\delta)= \Phi_I\cup\{\pm\delta\}$$
      since, from $L/Q=G/B$, every parabolic subgroup $R$ of $L$ that lies
      immediately over $Q$ has the property that $R/Q\cong\P^1$.
      So $\delta\in \Phi_I^\perp$, so that $\Phi_\Delta=\Phi_I\perp \Phi_{\Delta-I}$.
      So $I$ is empty and $Q$ is a Borel subgroup of $L$.
      Therefore
      $$\dim G =2\dim G/B +\rk(G)=2\dim L/Q+\rk(L)
      =\dim L$$
      (since the rank of a semi-simple group equals the Picard number of its flag variety).
      But $G\subseteq L$.
    \end{proof}
  \end{theorem}
  \begin{theorem} If $G$ and $H$ are adjoint group schemes
    over $\Sigma$ and $C(G)=C(H)$ then $G$ and $H$ are locally
    isomorphic over $\Sigma$.
    \begin{proof} $\sF^G$ and $\sF^H$ are locally isomorphic,
      by Theorem \ref{9.1},
      and $G=\Aut^0(\sF^G)$.
    \end{proof}
  \end{theorem}
  \begin{theorem}\label{conn} If $\Sigma$ is a geometric point,
    $g\in \Aut(\sF^G)$ and
    acts trivially on $\Pic(\sF^G)$ then $g\in G$.
    \begin{proof} Put $\sF^G=\sF$. The cohomology ring
      $H^*(\sF, \Q_{\ell})$ is generated as a $\Q_{\ell}$-algebra by
      $\Pic(\sF)$. So $g$ acts trivially on $H^*(\sF, \Q_{\ell})$
      and then (Lefschetz--Verdier) fixes a point of $\sF$. That
      is, we can regard $g$ as an automorphism of $G$ that preserves
      some Borel subgroup $B$. Choose
      a maximal torus $T$ in $B$; then $T$ is maximal in $G$
      and $g(T)$ is $B$-conjugate to $T$. So we can assume that $g$
      preserves $T$ and $B$. Since $\X^*(T)\subseteq\Pic(\sF)$ the action of
      $g$ on $T$ is trivial.
      Then for every root $\b$ the associated root subgroup
      $U_\b$ given by \cite{SGA3} XXII Th. 1.1 is preserved by $g$.
      Then we can replace $g$ by $t g$ for some
      $t\in T$ to make $\ad(g)$ act on $U_\a$
      as the identity for every simple $\a$.

      For any simple root $\a$ let $S_\a$
      denote the copy of
      $(P)SL_2$ in $G$ that is generated by $U_\a$ and $U_{-\a}$;
      then $g$ induces an automorphism $g_\a$ of $S_\a$
      which acts trivially on the
      subgroup of diagonal matrices and on $U_\a$.
      Every automorphism of $(P)SL_2$
      is inner, and then a calculation with $2\times 2$ matrices
      shows that $g_\a$ is the identity.
      The groups $S_\a$ generate $G$
      and so $g=1$.
    \end{proof}
  \end{theorem}
    We say that $G$ is \emph{quasi-split by $y$} if
    $y\in\sF^G(\Sigma)$ and that
    $G$ is \emph{pre-split by $y$} if it is quasi-split
    by $y$ and the {\'e}tale
    sheaf $M= \Pic^G(\sF^G)/\Pic(\Sigma)$ of finite
    free $\Z$-modules on $\Sigma$ is constant.
    Since $M$ is the character group of the
    constant torus $T\times\sF^G$ which is the
    reductive quotient of the universal Borel subgroup
    $\sB\to\sF^G$
    of $G\times\sF^G$ any
    torus in a Borel subgroup of a pre-split group is split.

    Assume that $G$ is adjoint
    (that is, $G$ acts effectively on $\sF^G$)
    and pre-split by $y\in\sF^G(\Sigma)$.
    Put $B=B_G=\Stab(y)$ and let $U=U_G$ denote the unipotent
    radical of $B$. 
    Put $X^G_S=\cup\ \tau_i^{-1}(\tau_i(y)),$
    the union of the Schubert curves through $y$.
    There is a surjective homomorphism
    $B\to \Aut^0 (X^G_S)$
    whose kernel is $[U, U]$.

    Suppose also that $H$ is an adjoint
    group over $\Sigma$, pre-split by
    $z\in\sF^H(\Sigma)$ and that $C(H)=C(G)$.
  
  \begin{theorem}\label{gamma}
    Assume that
    $\psi_S:X^G_S\to X^H_S$ is a stratified isomorphism
    and that $H^1(\Sigma, -\a)=0$
    for every $\a\in\Phi_{++}$.
    
    Then there is an isomorphism $\psi_\sF:\sF^G\to \sF^H$
    that extends $\psi_S$ 
    and an isomorphism $\psi:G\to H$ that induces $\psi_\sF$.
    Each set of such isomorphisms is a torsor under
    $H^0(\Sigma, [U, U])$.
    \begin{proof} The scheme
      $I=\Isom^C_\Sigma(\sF^G,\sF^H)$,
      where the superscript means ``isomorphisms that preserve $C$'',
      is a torsor
      under $G\to\Sigma$. As such $I$ is a class $\xi$
      in $H^1(\Sigma, G)$. The data of $y$ and $z$ lift $\xi$ to
      $\eta\in H^1(\Sigma, B)$. 
      The datum of $\psi_S$
      lifts $\eta$ to a class $\zeta\in H^1(\Sigma, [U, U])$;
      our assumptions ensure that
      this $H^1$ is trivial.
    \end{proof}
  \end{theorem}
  \begin{remark} The $2$-dimensional
    Schubert subscheme $\cup_{\ell(w)\le 2}X_w$
    determines the Cartan matrix, as explained in
    Sections \ref{Cartan} and \ref{six}, but is not enough
    to rigidify the situation because $B$
    does not always act effectively on it.
  \end{remark}
  
  Assume now also that $T_G\subset B_G$
  and $T_H\subset B_H$ are 
  tori which, everywhere locally on $\Sigma$,
  are maximal in $G$ and $H$, respectively.
  As already remarked, these tori are split, so that $G$ and $H$ are
  split over $\Sigma$. Conversely, any split group is also pre-split
  (for example, by a choice of point fixed under the torus action).
  
  \begin{theorem} Suppose that $\psi_S:X^G_S\to X^H_S$
    and $\psi_T:T_G\to T_H$ are isomorphisms which are compatible
    with the actions of $T_*$ on $X^*_S$ for $*=G, H$.
    Then there are unique
    isomorphisms $\psi_\sF:\sF^G\to\sF^H$ and $\psi:G\to H$
    that extend the datum $(\psi_S,\ \psi_T)$.
    \begin{proof}  In the notation
      of the proof of Theorem \ref{gamma}
      it is enough to show that $\zeta$ is trivial.
      The assumptions on the tori imply
      that $\zeta$ lifts to $\tau\in H^1(\Sigma, \G)$,
      where $\G$ is the subgroup of $[U_G, U_G]$ defined
      by the requirement that it should
      normalize $T_G$. Then
      $\G\subseteq [U_G, U_G]\cap N_G(T_G)=1$
      so that $\tau$, and therefore $\zeta$, is trivial.
    \end{proof}
  \end{theorem}
  Restating this in terms of groups over $\Z$
  gives the uniqueness theorem of \cite{SGA3}
  for adjoint groups.

  Fix a Cartan matrix $C$. There is an adjoint group
  $G_\Z$ that belongs to $C$ and is pre-split, say by
  $0\in\sF_\Z=\sF_{G_\Z}$. There is a split torus
  $T_\Z\subseteq B_\Z=\Stab(0)$ that is maximal in $G_\Z$,
  and so $G_\Z$ is split.

  On the other hand, take 
  the union $P$ of $r$ copies of $\P^1_\Z$
  that are identified at $0$ with relative
  embedding dimension $r$.
  Fix an isomorphism $h:P\to X^{G_\Z}_S$
  such that $h(0)=0$, the $i$'th copy of $\P^1_\Z$ maps
  to $\tau_i^{-1}(\tau_i(0))$ and the fixed points of $T_\Z$
  on each $\P^1_\Z$ are $0,\infty$. (The existence
  of $h$ follows from the fact that the curves
  $\tau_i^{-1}(\tau_i(0))$  have independent tangent directions
  at $y$, which in turn follows from the fact
  that the simple roots are linearly independent
  in $\Lie(B)/\Lie(U)$.)

  Suppose now that $G\to\Sigma$ is also adjoint, with
  Cartan matrix $C$, and pre-split by $y\in\sF^G(\Sigma)$.
  Suppose also that $T\subseteq B=\Stab(y)$ is a split torus,
  maximal in $G$ (so that $G$ is split), and that $j:P\times\Sigma \to X^G_S$
  is an isomorphism such that the $i$'th copy of
  $\P^1_\Sigma$ maps to $\tau_i^{-1}(\tau_i(y))$ and the fixed points of
  $T$ on each $\P^1_\Sigma$ are $0,\infty$. These data
  are a \emph{pinning} of the split group $G$ and there is a unique isomorphism
  $\psi_T:T_\Z\times\Sigma \to T$ that intertwines $j$
  with the actions of $T_\Z$ on $P$ and of $T$ on $X^G_S$.

  \begin{theorem}\label{unique} (Uniqueness for adjoint groups)
    There is a unique isomorphism
    $\psi:G_\Z\times\Sigma\to G$ such that $\psi$
    restricts to $\psi_T$, $\psi_\sF(0)=y$
    and $\psi_\sF\circ j=h\times 1_\Sigma$.
    \noproof
  \end{theorem}
  We recover the general uniqueness theorem
  as follows. Fix a pinned root datum (\cite{SGA3} XXIII 1.5)
  $\sR=(M, \Phi, \Delta, M^\vee, \Phi^\vee)$ and suppose that
  $G$ is a pinned reductive group over $\Sigma$ whose pinned root datum
  is $\sR$. Recall that then $M=\Pic^G(\sF^G)$ and $\Phi$ is the
  $W$-orbit of the simple roots, which are the relative tangent bundles of the
  $\P^1$-fibrations $\tau_i$.
  There is a pinned Chevalley group $G^\sR_\Z$
  over $\Sp\Z$
  whose pinned root datum is $\sR$.
  \begin{theorem}\label{uni} (Uniqueness in general)
    There is a unique isomorphism
    $G\to G^\sR_\Z\times_{\Sp \Z}\Sigma$ of pinned 
    groups.
    \begin{proof} We have proved the result for adjoint groups.
      We continue as follows.
      \begin{enumerate}
      \item From $G^{ad}$ we can construct its universal cover
        $G^{sc}$ as follows. Write $\sF=\sF^G=\sF^{G^{ad}}$.
        Let $T_1$ be the torus with $\X^*(T_1)=\Pic(\sF)$;
        then there is a universal torsor
        $\sT_1\to\sF$ under $T_1$. That is,
        a character $\chi$ of $T_1$ defines a line bundle
        $\sT_1\times^\chi\A^1\to\sF$, and this defines
        an isomorphism $\X^*(T_1)\to\Pic(\sF)$.
        Then there is
        a central extension 
        $$1\to T_1\to\sG\to G^{ad}\to 1$$
        where $\sG$ is defined as a group-valued functor by
        $$\sG=\{(\lambda, a) \vert  a\in G^{ad}\ \textrm{and}\
        \lambda : a^*\sT_1\buildrel\cong\over\to\sT_1\
        \textrm{is\ an\ isomorphism}\}$$
        and the group law is
        $(\lambda, a).(\mu, b)=(\mu\circ b^*\lambda, ab).$
        Since $\sG\to G^{ad}$ is relatively representable
        (the fibre over a point of $G^{ad}$ is represented by
        a torsor under $T_1$)
        $\sG$ is representable. Moreover, the scheme $\sG$
        is then a $T_1$-bundle over a smooth affine scheme, so is itself
        smooth and affine. 
        We then construct $G^{sc}$
        as the derived subgroup of $\sG$ and see that $G^{sc}$
        is characterized in terms of $G^{ad}$ by the fact that 
        it is a central extension of $G^{ad}$ by a subgroup of a torus
        and every line bundle on $\sF$
        is uniquely $G^{sc}$-linearized.

      \item $\sR^{ad}= (Q, \Phi, \Delta, Q^\vee, \Phi^\vee)$
        and so determines the simply connected
        root datum $\sR^{sc}$ from the formula
        $\sR^{sc}=(Q^\vee, \Phi, \Delta, Q, \Phi^\vee).$

      \item An arbitrary semi-simple
        $G$ is trapped between $G^{sc}$ and $G^{ad}$
        and is determined by the group $\Pic^G(\sF)$.
        This group is just $M$ and so is determined by $\sR$.
        So the result is proved for semi-simple groups.

      \item\label{ss to red} Let $Z_1$ be the torus subgroup scheme of the centre
        $Z$ of $G$, defined by the      property
        that $\X^*(Z_1)$ equals $\X^*(Z)$ modulo its
        torsion subgroup, and $G^{der}$
        the derived subgroup of $G$. Then $G^{der}$ is semi-simple
        and there
        is a canonical central isogeny $\pi:G^{der}\times Z_1\to G$.
        Since $\ker\pi$ is determined by $\sR$
        we are done.
      \end{enumerate}
    \end{proof}
  \end{theorem}
\end{section}
\begin{section}{The homomorphism theorem}
  Suppose that $G$ and $H$ are pinned reductive
    groups over $\Sigma$ with pinned root data
    $\sR(G)=(M, \Phi,\Delta, M^\vee, \Phi^\vee)$ and
    $\sR(H)=(M', \Phi', \Delta', M'^\vee,  \Phi'^\vee)$.
    \begin{theorem} Each
    morphism $\phi:\sR(G)\to\sR(H)$
    of pinned root data
    is induced by a unique homomorphism $G\to H$.
    \begin{proof} Comprised in $\phi$ is a $\Z$-linear
      map $f:M'\to M$. Define multiplicative groups $X$ and $Y$
      by $\X^*(X)=\coker f$ and $\X^*(Y)=\ker f$.
      Then $X$ is central in $G$, there is a surjection
      $H\to Y$ and the groups $G/X$ and $\ker (H\to Y)$
      have isomorphic pinned root data.
      Theorem \ref{uni} provides a unique isomorphism
      $G/X\buildrel\cong\over\to \ker(H\to Y)$
      compatible with the isomorphism of pinned root data
      and completes the proof.
    \end{proof}
  \end{theorem}
  \begin{theorem} Suppose that $\Sigma$ is of characteristic
    $p>0$. Then each $p$-morphism $\phi:\sR(G)\to\sR(H)$ of
    pinned root data is induced by
    a unique homomorphism
    $G\to H$.
    \begin{proof} 
        We shall understand a $p$-morphism
      of pinned root data to comprise
      a $\Z$-linear map
      $f:M'\to M$,
      a bijection $u:\Delta\buildrel\cong\over\to \Delta'$
      and a map
      $q:\Delta\to\{p^n : n\in\N\}$
      such that $f(u(\a))=q(\a)\a$
      and ${}^tf(\a^\vee)=q(\a)u(\a)^\vee$ for all $\a\in \Delta$. In particular,
      $\phi$ induces an isomorphism of Coxeter systems,
      so that $f$ is $W$-equivariant and
      $u, q$ extend uniquely from $\Delta$ to $\Phi$,
      as required in the standard definition of a $p$-morphism
      \cite{SGA3} XXI 6.8.
      Because $q(\a)(\a.\b^\vee)=q(\b)(u(\a).u(\b)^\vee)$
      one Cartan matrix determines the other.
      
      Suppose first that $G$ and $H$ are both adjoint.
      We shall give two proofs in this case. The first
      uses the Bott--Samelson $\Sigma$-schemes
      $Z^G_{\bw}$ and $Z^H_{u({\bw})}$
      and the    automorphism group schemes $\sK^G_{\bw}$
    and $\sK^H_{u({\bw})}$ of \ref{SectPrelim},
    where ${\bw}=s_{i_1}\ldots s_{i_n}$
    and $s_j=s_{\a_j}$. The second proof goes via
    the construction of a certain normal
    subgroup $K$ of $G$ and then proving that
    $G/K$ and $H$ are uniquely isomorphic.
    
    Here is the first proof, for adjoint groups.
    
    Consider the fundamental dominant weights $\varpi_i^G$
    and $\varpi_i^H$;
    they satisfy $(\varpi^G_i.\a_j^\vee)=(\varpi_i^H.u(\a_j)^\vee)=\delta_{i j}$.
    Then $u(\a)=\sum(u(\a).u(\a_i)^\vee)\varpi_i^H$.
    
    Define $\chi^G_{{\bw}, a}\in\Pic(Z^G_{\bw})$ to be the pull back
    of $\varpi^G_{i_a}$ under the $a$th projection
    $Z^G_{\bw}\to\sF^G$. Then 
    the pull back
    of $\varpi^G_b$ to $Z^G_{\bw}$ under the final projection
    is $\chi^G_{{\bw}, j(b)}$, where $j(b)$ is the maximal
    integer $m$ such that $i_m=b$,
    if $s_b\in{\bw}$. If $s_b\not\in{\bw}$ then
    we define $\chi^G_{{\bw}, j(b)}$ to be trivial.
    So if $\lambda=\sum n_b\varpi^G_b$ then
    $$\lambda\vert_{Z^G_{\bw}}=\sum_{s_b\in{\bw}}n_b\chi^G_{{\bw}, j(b)}.$$
    
    We shall construct, by induction on $\ell({\bw})$,
    a radicial morphism $g_{\bw}:Z^G_{\bw}\to Z^H_{u({\bw})}$
    such that
    $$g_{\bw}^*\chi^H_{u({\bw}), a}=q(\a_{i_a})\chi^G_{{\bw}, a}$$
    and $g_{\bw}$ is equivariant under
    a homomorphism $\sK^G_{\bw}\to \sK^H_{u({\bw})}$.

    Suppose that ${\bw}={\bv}s_\a$ (so $\a=\a_{i_n}$)
    and that $g_{\bv}$ has
    been constructed with the required properties.
    Since $q(\a)(\a.\b^\vee)=q(\b)(u(\a).u(\beta)^\vee)$
    we get, by a straightforward substitution,
    $$g_{\bv}^*(u(\a)\vert_{Z^H_{\bv}})=q(\a)\a\vert_{Z^G_{\bv}}.$$
    
    The fixed locus of $B^H$ on $H^1(Z^H_{u({\bv})},-\a)$
    is contained in the fixed locus of $\sK_{u({\bv})}^H$,
    and so pulls back under $g_{\bv}$ to a subspace of
    $H^1(Z^G_{\bv}(v),-q(\a)\a)^{B^G}$,
    by the inductive equivariance, which takes $B^G$ to
    a subgroup of $\sK_{u({\bv})}^H$.
    Therefore there is a commutative diagram
    with Cartesian square ($\Pi$ is the fibre product)
    $$
    \xymatrix
    {
      {Z^G_{\bw}}\ar[dr]\ar@/^1pc/[rr]^{g_{\bw}}\ar[r]_{F}&
      {\Pi}\ar[r]\ar[d]
      \ar@{}|-{\square}[dr]&
      {Z^H_{\bw}}\ar[d]\\
      &{Z^G_{\bv}}\ar[r]_{g_{\bv}}&
      {Z^H_{\bv}}
    }
    $$
    where $F$ is the relative $q(\a)$ Frobenius;
    this defines $g_{\bw}$.

    We must verify that
    $g_{\bw}^*\chi_{u({\bw}), a}^H=q(\a_{i_a})\chi_{{\bw}, a}^G$.
    It suffices to compare them on the marked section
    $Z^G_{\bv}\inj Z^G_{\bw}$ and on a geometric fibre $\a^\vee$
    of $Z^G_{\bw}\to Z^G_{\bv}$. If $a<n$ both
    divisor classes have the same restriction to $Z^G_{\bv}$,
    by induction, and both are trivial on $\a^\vee$.

    If $a=n$ then they agree on $\a^\vee$,
    by construction. Moreover,
    $g_{\bv}^*\chi^H_{u({\bw}), i_n}=g_{\bv}^*\chi^H_{u({\bv}), j'(i_n)},$
    which, by induction, is $q(\a)\chi^G_{{\bv}, j'(i_n)}=
    q(\a)\chi^G_{{\bw}, i_n}\vert_{Z^G_{\bv}}.$
    Here $j'(i_n)$ is defined for $\bv$ as $j(i_n)$
    is for $\bw$.



      Assume that
      $\ell({\bw})\ge 1$ and write      ${\bw}=s_\b{\bx}$.
      Then $Z^G_{\bw} = (q^G_{\bx})^{-1}(\b^\vee)$,
      where $\b^\vee$ is a copy of $\P^1$ (a fibre
      of the $\P^1$-fibration $\sF^G\to\sP^G_\b$).
      That is, $Z^G_{\bx}$ is a fibre of a $G$-equivariant map
      $\tZ^G_{\bx}\to\sF^G$ and $Z^G_{\bw}$ is the inverse
      image of the curve $\b^\vee$ on which $P_\b$
      acts; therefore
      there is a homomorphism from the
      minimal parabolic subgroup $P^G_\b$ of $G$ to
      $\Aut^0(Z^G_{\bw})$
      and a commutative triangle
      $$
      \xymatrix
      {
        {P^G_\b}\ar[r]^-{r_{\bw}}\ar[dr]& {\Aut^0(Z^G_{\bw})}\ar[d]\\
        {}&{PGL_{2, \b}}
      }
      $$
      where $PGL_{2, \b}$ is the copy of $PGL_2$ that acts on
      $\b^\vee$. There is an analogous picture on the $H$ side
      and the homomorphism $\sK^G_{\bw}\to \sK^H_{u({\bw})}$
      covers an isogeny $PGL_{2, \b}\to PGL_{2, u(\b)}$.

      Observe the following things.
      \begin{enumerate}
        \item From its construction, the morphism
          $g_{\bw}$ is equivariant
          under $\Aut^0(Z^G_{\bw})$
          with respect to a homomorphism
          $t_{\bw}: \Aut^0(Z^G_{\bw})\to \Aut^0(Z^H_{u({\bw})})$
          that covers an isogeny $PGL_{2, \b}\to PGL_{2, u(\b)}$
          and not merely under $\sK^G_{\bw}$.
        \item If $\bw$ represents $w\in W$ there is a commutative diagram
          $$
          \xymatrix
          {
            {Z^G_{\bw}}\ar[r]^-{g_{\bw}}\ar[d]_{a^G_{\bw}}&
            {Z^H_{u({\bw})}}\ar[d]^{a^H_{u({\bw})}}\\
            {\sF^G_w}\ar[r]_-{\g_w}&
            {\sF^H_{u(w)}}
          }
          $$
          where $\g_w$ depends only on $w$.
          The reason is that, as in the proof of Proposition
          \ref{5.1}, the curves collapsed
          by $a^G_{\bw}$ are exactly the curves
          orthogonal to the divisor classes represented
          by the simple roots and the same is true of
          $a^H_{u({\bw})}$, so that $\g_w$ takes any curve
          collapsed by $a^G_{\bw}$
          to one that is collapsed by $a^H_{u({\bw})}$.
          Moreover, $\g_w$ is radicial since $g_{\bw}$ is radicial
          and the vertical maps are birational.
        \item If $\bw$ is reduced then $a^G_{\bw}$ and
          $a^H_{u({\bw})}$ are proper birational morphisms of normal
          varieties. Since also
          $P^G_\b$ is smooth and connected, this
          square is equivariant
          with respect to the homomorphism
          $t_{\bw}\circ r_{\bw}:P^G_\b\to \Aut^0(Z^H_{u({\bw})})$,
          as follows from the next lemma.
        \end{enumerate}
        \begin{lemma} If $g:X\to Y$ is a proper birational
          morphism of normal quasi-projective varieties
          and $P$ is a smooth connected 
          group acting on $X$, then $P$ acts
          on $Y$ and $f$ is $P$-equivariant.
          \begin{proof} It is enough to prove that $P$
            preserves the exceptional locus $E$ of $f$.
            As in the proof of Proposition \ref{5.1},
            this locus is covered by curves $C$ such that
            $C.D=0$ for every divisor class $D$ that pulls
            back from $Y$. Since $P$ is connected it preserves
            the numerical class of each such $D$, and
            so preserves $E$.
          \end{proof}
        \end{lemma}
        
        Suppose that ${\bw}_0$ begins with $s_\b$ and
        is a reduced expression for the
        longest element $w_0$ of $W$.
        Then $\g_{w_0}:\sF^G\to\sF^H$ is
        $P^G_{\b}$-equivariant with respect to a non-trivial
        homomorphism $\psi_\b:P^G_{\b}\to\Aut^0(\sF^H)=H$.
        But for every simple
        $\b$ there is a reduced expression ${\bw}_0$
        that begins with $s_\b$, and the groups
        $P^G_\b$ as $\b$ runs over the simple roots
        generate $G$.

        Since, as remarked above, $\g_{w_0}$ is radicial,
        there is, for some $n$, a factorization
        $$
        \xymatrix
        {
          {\sF^G}\ar[r]_{\g_{w_0}}\ar@/^1pc/[rr]^{F^n}&
          {\sF^H}\ar[r]&
          {\sF^G}
        }
        $$
        of $F^n$, the
        $n$th power of the Frobenius. Certainly
        $F^n$ is $G$-equivariant, and every $P^G_{\b}$
        preserves this factorization. Therefore
        $G$ preserves it, so that $\g_{w_0}$
        is equivariant with respect to some homomorphism
        $\psi:G\to H$. Then $G/\ker\psi\to H$
        is an isomorphism, and is uniquely determined
        by the pinned root data.

        This concludes the first proof for adjoint groups.
        \bigskip

        We now give the second proof for adjoint groups.
     
        We
        can take both root data to be irreducible.
      Define
      $\phi$ to be \emph{primitive} if $q(\a)=1$ for some
      $\a\in \Phi$; then $\phi$ is the composite
      of a primitive $p$-morphism and a constant one
      (one which is multiplication by a fixed
      power of $p$)
      (\cite{SGA3} XXI 6.8.6). If $\phi$ is primitive
      but not constant
      then $q(\a)\lng(\a)$ is constant,
      $q$ takes the values $\{1, p\}$,
      there are roots of two different lengths and the
      ratio of the different lengths is $p$ (\cite{SGA3} XXI 7.5.2).
      Note that  on
      the group side a constant $p$-morphism corresponds
      to a power of the Frobenius relative to $\Sigma$ so we can assume that
      $\phi$ is primitive but not constant. Then we construct
      $G\to H$ by first constructing its kernel $K$
      as a subgroup of the kernel $G_{(1)}$ of
      Frobenius acting on $G$. The action of $G$ on $G_{(1)}$
      by conjugation is equivalent to its adjoint action
      on $\Lie(G)$, so that a subgroup $K$ of $G_{(1)}$
      is normal in $G_{(1)}$ if and only if $\Lie(K)$ is a
      $p$-closed ideal of the $p$-Lie algebra
      $\Lie(G)$, and also $K$ is normal in $G_{(1)}$
      if and only if it is normal in $G$.
      
      Recall that $\Lie(G)=\mathfrak g_\Z\otimes\sO_\Sigma$
      where $\mathfrak g_\Z$ is Chevalley's Lie algebra over $\Z$
      \cite{Ch}. Define $\mathfrak K\subseteq \Lie(G)$ to be
      the $\sO_\Sigma$-span of $\Lie(T)$
      and the generators
      $e_\a$ for the short roots $\a$. 

      \begin{lemma} $\mathfrak K$
        is a $p$-closed ideal of $\Lie(G)$.
        \begin{proof} 
          By \cite{Ch}, p. 24,
        $[e_\a, e_\b]=\pm m_{\a\b} e_{\b+\a}$
        when $\a,\b$ and $\b+\a$ are roots
        and $m_{\a\b}$ is the smallest integer $m>0$ such that
        $\b-m\a$ is not a root.
        We can assume that $\a$ is short and
        $\b+\a$ is long. According to \cite{St}
        Lemma 3,
        if $I_\a(\b)=[\b-r\a,\b+s\a]$ then
        $r+1=s\lng(\a+\b)/\lng(\a)$, so that
        $m_{\a\b}=r+1=sp=0$ and $\mathfrak K$
        is an ideal.

        $\mathfrak K$ is $p$-closed because,
        if $\a$ and $\b$ are short roots, then
        $(\ad(e_\a))^2(e_\b)$ is zero or a multiple of $e_{2\a+\b}$
        and $2\a+\b$ is not a short root, so that
        $(\ad(e_\a))^2=0$ and then $(\ad(e_\a))^p=0$.
        $\Lie(T)$ is generated by $p$-idempotent
        elements, so the lemma is proved.
      \end{proof}
    \end{lemma}
    Take $K$ to be the subgroup of $G_{(1)}$,
    normal in $G$,
    with $\Lie(K)=\mathfrak K$. 
    \begin{lemma}\label{ifinj} $G/K$ is adjoint and is
      uniquely isomorphic to
      $H$ as a pinned group.
      \begin{proof} 
        The action of $K$ on $\sF^G$ is not free
        but, because $G$ acts transitively on $\sF^G$
        and $K$ is finite and normal in $G$,
        the stabilizer subgroup scheme
        of $K\times_\Sigma\sF^G$ is finite and flat
        over $\sF^G$. Therefore there is a
        quotient $\rho:\sF^G\to X=K\backslash\sF^G$
        (\cite{SGA3} V) that is smooth and projective over $\Sigma$
        and has a transitive action of $G/K$. Let $x\in\sF^G$
        be the point stabilized by the Borel subgroup
        $B$ of $G$ that is generated by $T$ and the $U_{-\a}$
        for the positive roots $\a$. Then the stabilizer
        of $\rho(x)$ is $B/B\cap K$ so is smooth, connected
        and soluble. Therefore $X$ is the flag variety of $G/K$.
        Moreover, the action of $G/K$ is effective, and so,
        by Lemma \ref{adj},
        $G/K$ is adjoint.

        For a simple root $\a$ let $P_\a\supset B$ be the corresponding
        minimal parabolic subgroup. Then, from the construction
        of $K$, $P_\a/B$ maps with degree $q(\a)$ to its
        image in $\sF^{G/K}$. Since the Frobenius $\sF^G\to\sF^G$ factors
        through $\sF^{G/K}$ the $\P^1$-fibrations of $\sF^G$ and of
        $\sF^{G/K}$ match up under $\rho$. This gives a bijection $v$
        from the set of simple roots in $\sF^G$ (the relative tangent
        bundles of the $\P^1$-fibrations) to those in $\sF^{G/K}$
        such that, if $\G_\a$ is a fibre (a simple coroot $\a^\vee$) then
        $\rho_*\G_\a=q(\a)\G_{v(\a)}.$
        There is a commutative diagram with Cartesian square
        ($\Pi$ is the fibre product)
        $$
        \xymatrix
        {
          {\sF^G}\ar[r]_\s\ar[dr]\ar@/^1pc/[rr]^{\rho}&
          {\Pi}\ar[r]\ar[d]\ar@{}|-{\square}[dr]&
          {\sF^{G/K}}\ar[d]\\
          &{\sP_\a}\ar[r]&{\sP_{v(\a)}}
        }
        $$
        where
        either $\s$ is an isomorphism, in which case
        $\rho^*v(\a)=\alpha$, or $\s$ is the relative
        Frobenius, in which case $\rho^*v(\a)=p\a$.
        So $\rho^*v(\a)=\lambda(\a)\a$
        for some $\lambda(\a)\in\N$.        
        In either case the projection formula
        $(\rho^*v(\a).\G_\b)=(v(\a).\rho_*\G_\b)$
        gives, by taking $\a=\b$, $\lambda(\a)=q(\a)$
        and then in general
        $$q(\a)(\a.\b^\vee)=q(\b)(v(\a).v(\beta)^\vee).$$
        It follows that $H$ and $G/K$ are pinned adjoint groups
        with isomorphic pinned root data and
        so are uniquely isomorphic.
      \end{proof}
    \end{lemma}
    This concludes the second proof
    for adjoint groups.
        
    Next, suppose that $G$ and $H$ are semi-simple.
    We have a homomorphism $G^{ad}\to H^{ad}$, and so
    homomorphisms $G^{sc}\to H^{sc}\to H$. The kernel of
    $G^{sc}\to H$ consists of those $g\in G$ such that
    $g$ acts trivially on $\sF^G$ and on every
    $\rho^*\sL$  where $\sL$ is an $H$-linearized
    line bundle on $\sF^H$. This is just the multiplicative
    group $S$ defined by $\X^*(S)=\coker(M'\to N)$,
    where
    $N=\Pic^{G^{sc}}(\sF^G)\supseteq M=\Pic^G(\sF^G)$
    and so $G^{sc}\to H$ factors through $G$
    to give a homomorphism $G\to H$.
    The induced homomorphism $G^{ad}\to H^{ad}$
    is unique, so $G\to H$ is unique.
    
    The theorem is now proved for semi-simple groups.

    To extend to reductive groups in general we
    argue as in the proof of Theorem \ref{uni} \ref{ss to red}.
  \end{proof}
\end{theorem}
\end{section}

\end{document}